\magnification=1100
\baselineskip=15pt

\hsize=150mm
\vsize=205mm

\let \al=\alpha
\let \be=\beta
\let \var=\varphi
\let \vare=\varepsilon

\let \de=\delta

\let \la=\lambda

\let \p=\partial
\let \q=\quad
\let \qq=\qquad
\let \med=\medskip
\let \smal=\smallskip

\def\R{{\rm I\kern
-1.6pt{\rm R}}}
\def\C{{\rm |\kern
-4.6pt{\rm C}}}
\def\N{{\rm I\kern
-4.0pt{\rm N}}}

\def\eq{\eqno}
\def\q{\quad}

\def\ter{\hfill \vrule width 5 pt height 7 pt depth - 2 pt\smallskip}

\def\system#1{\left\{\null\,\vcenter{
\ialign{\strut\hfil$##$&$##$\hfil&&\enspace$##$\enspace&
\hfil$##$&$##$\hfil\crcr#1\crcr}}\right.}

\

\centerline{\bf Global Attractivity and Extinction for  Lotka-Volterra systems}

\centerline{\bf  with infinite delay and feedback controls}

\

\vskip0.3in

\centerline{
Teresa Faria \parindent 0cm 
{\footnote{$^{\star}$}{Corresponding author.Fax:+351 21 795 4288; Tel: +351 21 790 4929.}}$^,$
{\footnote {$^{\rm a}$} {The research was supported by Funda\c c\~ao para a Ci\^encia e a Tecnologia (Portugal),  PEst-OE/\-MAT/\-UI0209/2011.} }
and Yoshiaki Muroya \parindent 0cm{\footnote{$^{\rm b}$}{The research was supported by Scientific Research (c), No.24540219 of Japan Society for the Promotion of Science.} }
}

\vskip .2cm
\centerline{$^{\rm a}$Departamento de Matem\'atica  and CMAF, Faculdade de
Ci\^encias,}
\centerline{Universidade de Lisboa, Campo Grande, 1749-016 Lisboa,
Portugal}
\centerline{tfaria@ptmat.fc.ul.pt}

\vskip0.3in

\centerline{$^{\rm b}$Department of Mathematics, Waseda University}
\centerline{ 3-4-1 Ohkubo, Shinjuku-ku, Tokyo 169-8555, Japan}

\centerline{ymuroya@waseda.jp}

\vskip0.3in

{\bf Abstract}: The paper deals with a
multiple species Lotka-Volterra model with infinite distributed delays and feedback controls, for which we assume a weak form of diagonal dominance of the instantaneous  negative intra-specific terms over the infinite delay effect in both the population variables and  controls. General sufficient conditions for the existence and attractivity of  a saturated equilibrium are established. When the saturated equilibrium is on the boundary of $\R^n_+$, sharper criteria for the extinction of all or part of the populations are given.
While  the literature usually   treats the case of competitive systems only, here no restrictions on the signs of the intra- and inter-specific delayed terms are imposed. Moreover, our technique does not require the construction of Lyapunov functionals.

\vskip0.3in


\vskip0.3in


\vfill\eject

{\bf 1. Introduction}

\med

After several decades of intensive  study and use of functional differential equations (FDEs) in population dynamics, it is now very well understood that   the introduction of delays in differential equations  leads in general to more realistic population models,  and  much more complex and rich dynamics.
Nevertheless, delays are not
 harmless and  often
 create instability and oscillations,  unless they  are  either small  or neutralized by instantaneous terms. When the delays are  infinite, it is not clear how to surpass the  effect of the infinite  past of the system, so in order to obtain stability results some form of  instantaneous dominance is expected.
On the other hand, the consideration of  FDEs with infinite delay is relevant to account for  systems with ``infinite memory", and  goes back to the works of Volterra.
In fact, for 
Lotka-Volterra systems or other general population models, whether the global stability may persist under large or even infinite delays without strictly dominating instantaneous negative feedbacks is a question that has attracted the interest of many researchers, and  had partial positive answers, see e.g. Kuang [15], Xu et al.~[26],  Faria [4], also for further references.

Recently, the study of population models with delays and controls, in particular Lotka-Volterra models, 
 has received some attention, see e.g.  [3,~8,~16,~20,~21,~24,~27], and references therein.
In this paper, we consider the following n-species Lotka-Volterra system with feedback controls and infinite delays:
$$\system{
\displaystyle x_i^{\prime}(t)&=x_i(t)\biggl(b_i-\mu_ix_i(t)-\sum_{j=1}^n a_{ij} \int_0^{\infty}K_{ij}(s)x_j(t-s)\, ds-c_i \int_0^{\infty}G_i(s)u_i(t-s)\, ds\biggr), \cr
\displaystyle u_i^{\prime}(t)&=-e_i u_i(t)+d_i x_i(t), \quad i=1,2,\ldots,n , \cr}
\eq(1.1)$$
where: $ \mu_i,  c_i , d_i, e_i$ are positive constants, $b_i,a_{ij}\in\R$, and
the kernels $K_{ij}, G_i: [0,\infty) \to [0,\infty)$ are $L^1$ functions, normalized so that 
$$
\int_0^{\infty} K_{ij}(s)\, ds=1,\q \int_0^\infty G_i(s)\, ds=1,
$$
for    $i,j=1,2, \ldots,n$.
Without loss of generality, we assume that for all $i$  the linear operators defined by $L_{ii}(\var)=\int_0^{\infty} K_{ii}(s)\var (-s)\, ds$, for bounded continuous  functions $\var:(-\infty,0]\to\R$, are non-atomic at zero, which amounts to have $K_{ii}(0)=K_{ii}(0^+)$.

In biological terms,   $x_i(t)$ denotes the density of the population  $i$ with Malthusian growth rate $b_i$ and instantaneous self-limitation coefficient $\mu_i>0$, and $a_{ii}$ and $a_{ij}(i\ne j)$ are respectively the intra- and inter-specific delayed acting coefficients; $u_i(t)$ denotes a feedback control variable, $i,j=1,2, \ldots,n$.
Due to the biological interpretation of model (1.1), we are only interested in positive  (or non-negative) solutions. We therefore consider solutions of system (1.1) with {\it admissible} initial conditions, i.e., 
$$
x_i(\theta)=\varphi_i(\theta)\ge 0, u_i(\theta)=\psi_i(\theta)\ge 0, \ \theta \in (-\infty,0), \quad \varphi_i(0)>0,\psi_i(0)>0,
\eq(1.2)$$
with $\varphi_i,\psi_i$   {\it bounded} continuous functions on $(-\infty,0],\  i=1,2,\ldots,n$. 

In order to have an effective   feedback control, it is natural to impose  that each $G_i\ (1\le i\le n)$ satisfies 
$$\int_0^{\infty} G_i(s) \var(-s)\, ds>0,
\eq(1.3)$$
for any positive, bounded,  continuous function $\var$ defined on $(-\infty, 0]$ with $\var(0)>0$. In particular, (1.3) holds if $G_i$ is continuous at zero with $G_i(0)>0$, or if
$G_i$ has a jump discontinuity at 0 with $G_i(0)-G_i(0^+)>0$.

  For simplicity of exposition, we consider systems (1.1), but our study applies to more general systems of the form
$$\system{
\displaystyle x_i^{\prime}(t)&=x_i(t)\biggl(b_i-\mu_ix_i(t)-\sum_{j=1}^n a_{ij} \int_0^{\infty}x_j(t-s)\, d\eta_{ij}(s)-c_i \int_0^{\infty}u_i(t-s)\, d\nu_i(s)\biggr), \cr
\displaystyle u_i^{\prime}(t)&=-e_i u_i(t)+d_i x_i(t), \quad i=1,2,\ldots,n , \cr}
\eq(1.4)$$
where all the coefficients  are as in (1.1), and
$\eta_{ij},\nu_i:[0,\infty)\to \R$ are  bounded variation  functions which are supposed to be normalized so that  their total variation is 1, and $\nu_i$ are  non-decreasing on $[0,\infty)$. Note that in (1.1) we supposed that $K_{ij}(t)\ge 0$ on $[0,\infty)$, but the above  scenario does not impose this restriction. 
Some of our general results  however require that the kernels $K_{ij}$ in (1.1) are non-negative,
or that the functions $\eta_{ij}$ in (1.4) are non-decreasing, although they  can be easily adapted to deal with systems  without such constraints.


Our study was strongly motivated by some previous  works of the present authors. The uncontrolled Lotka-Volterra system with infinite distributed delays was studied by Faria [4], and questions of partial survival and extinction of species in non-autonomous delayed Lotka-Volterra systems were addressed by Muroya in [19], see also [18]. The works of
Gopalsamy and Weng [8] and Li et al. [16], where special cases of  two-dimensional competitive  Lotka-Volterra systems with controls and no diagonal delays were studied, were an important source of inspiration for the present paper.  Here, the investigation refers to  controlled Lotka-Volterra models of any dimension $n$.  While  the literature usually  only deals with  the case of competitive systems  (i.e., systems with $a_{ij}\ge 0$ for $j\ne i$) with $b_i>0$, here no restrictions on the signs of $a_{ij}$ and $b_i$ will be imposed. Moreover,
 infinite delays are incorporated in the controls terms,  see also [21] for a competitive model. 
Another novelty is that our method does not require the construction of a specific Lyapunov functional.

Clearly, the introduction of controls in a delayed Lotka-Volterra system might change the existence, position, and stability  of equilibria. The main goal of the present paper is to address the global asymptotic dynamics of solutions to (1.1)-(1.2), in what concerns  establishing sufficient conditions for  the existence and  attractivity of a {\it saturated equilibrium} (see [14,~15] and Section 3 for a definition) not necessarily positive. 
As in  previous works [4,~8,~16,~21], we assume that system (1.1) satisfies some  form of diagonal dominance of the instantaneous  negative terms $\mu_ix_i(t)$ over the infinite delay terms, involving both the population variables  and the controls, so that the usual instability caused by the introduction of the delays is canceled.  For some of  our stability results, another prerequisite
   is that the uncontrolled Lotka-Volterra, system (1.1) with $c_i=0 \, (1\le i\le n)$, possesses already a globally attractive saturated equilibrium. 
These assumptions, although they seem restrictive,   are quite natural; moreover, here the main goal is to use the controls  to change the position of the saturated  equilibrium keeping its stability, as emphasised by some examples.  For a biological interpretation of the use of controls, see e.g. [8,~24,~27], also for additional  references.

   

\med

We now describe briefly  the contents of the paper.
From a theoretical  perspective, dealing with FDEs with infinite delays requires  a careful choice of a suitable  Banach phase space (usually called a {\it fading memory space}), in order to recover classical results of well-posedness of the initial value problem, existence and uniqueness of solutions, continuation of solutions, etc. 
For this reason, in   Section 2 we set some basic notation  for FDEs with infinite delays, and insert system (1.1) into such a framework. 
In Section 3, after studying the existence of a unique  saturated equilibrium  $(x^*,u^*)$ and the boundedness of positive solutions to (1.1), Theorems 3.2 provides a general criterion for the global attractivity of $(x^*,u^*)$. Also, a sufficient condition for the dissipativeness of (1.1) is given. In  Section 4,    sharper criteria are established  for the global attractivity  of a  saturated equilibrium  $(x^*,u^*)$  which is not strictly positive. In this situation, 
 this means  the extinction of all or part of the populations.    Our results turn out to be particularly powerful for predator-prey models.   We also emphasize that, for the uncontrolled system, we derive better results for partial (or total) extinction than the ones in [4]. Our techniques also allow to obtain a perturbation result for non-autonomous Lotka-Volterra systems with a limiting model of the form (1.1) or (1.4), as $t\to\infty$.   Section 5 is devoted to the particular case of  a 2-dimensional Lotka-Volterra system.  In Sections 4 and 5, some examples illustrate our results.

\bigskip

{\bf 2. An abstract formulation}

 \med
 
Since system (1.1) has unbounded delays, we must carefully formulate the problem by defining an appropriate Banach phase space where the problem is well-posed. 

 Let $g$ be a
function satisfying the following properties:
\med

{\parindent 1.5cm

\item {(g1)} $g:(-\infty ,0]\to [1,\infty)$ is a non-increasing continuous function and $g(0)=1$;

\item {(g2)} $\displaystyle{\lim_{u\to 0^-}{{g(s+u)}\over {g(s)}}=1}$ uniformly on $(-\infty ,0]$;

\item{(g3)} $g(s)\to \infty$ as $s\to -\infty$.}

\med
\noindent 
For $n\in\N$, define  the Banach space $UC_g=UC_g(\R^n)$,
$$\
UC_g:=\Big\{ \phi\in C((-\infty, 0];\R^{n}) : \sup_{s\le 0}{{|\phi(s)|}\over {g(s)}}<\infty,
{{\phi(s)}\over {g(s)}}\ {\rm is\ uniformly\ continuous\ on}\ (-\infty, 0]\Big\},$$
with the norm
$$\|\phi\|_g=\sup_{s\le 0}{{|\phi(s)|}\over {g(s)}},$$
where $|\cdot|$ is a chosen norm  in $\R^n$. 
Consider also  the space $BC=BC(\R^n)$ of bounded  continuous functions $\phi:(-\infty, 0]\to
\R^n$.
It is clear that $BC\subset UC_g$, with $\|\phi\|_g\le \|\phi\|_\infty$ for $\phi\in BC$ and  $\|\cdot\|_\infty$  the supremum norm in $BC$.
Here, $BC$ will be considered as
a subspace of $UC_g$, so $BC$ is endowed with the norm of $UC_g$.

The space $UC_g$ is an admissible phase space for $n$-dimensional FDEs with infinite delay
(cf. [12,~13]) written in the abstract form
$$\dot x(t)=f(t,x_t),\eq(2.1)$$
where $f:D\subset \R\times UC_g\to \R^n$ is continuous and, 
as usual, segments of solutions in the phase space $UC_g$ are denoted by $x_t$, $x_t(s)=x(t+s), s\le 0$, with components $x_{t,i}$.
Therefore, the standard results on existence and uniqueness of solutions for the Cauchy problem  $\dot x(t)=f(t,x_t), x_0=\var$ hold  when $f$ is regular enough and $\var\in BC$. Moreover, since $UC_g$ is a fading memory space, bounded positive orbits are precompact in $UC_g$ [9].

\med

We now set an appropriate formulation for problem (1.1)-(1.2).
From [10] and [6, Lemma 4.1], for any $\delta>0$ there is a continuous function $g$ satisfying (g1)--(g3) and such that
$$
\int_0^{\infty} g(-s)K_{ij}(s)\, ds<1+\de,\ \int_0^{\infty} g(-s)G_i(s)\, ds<1+\de, \q i,j=1,\dots,n.
\eq(2.2)$$
When dealing with systems (1.4), where the more general linearities  are  given by bounded variation functions $\eta_{ij}(s),\nu_i(s)$ with total variation 1 and $\nu_i(s)$ non-decreasing,  the above formulas (2.2) should be replaced by
$$\int_0^{\infty} g(-s)\, d|\eta_{ij}(s)|<1+\de,\ \int_0^{\infty} g(-s)\, d\nu_i(s)<1+\de,\q i,j=1,\dots,n. 
\eq(2.3)$$

Whenever an abstract setting is required, in what follows we shall always assume that (1.1) takes the abstract form (2.1) in the phase space $UC_g=UC_g(\R^{2n})$, for some fixed $\de >0$ and function $g$ satisfying (g1)--(g3) and (2.2), and consider solutions  with initial conditions 
$$
x_0=\var , \q u(0)=\psi,
\eq (2.4)$$
 where $(\var,\psi)\in BC(\R^{2n})$.
  System (1.1) has a unique solution $(x(t), u(t))=(x(t;\var,\psi), u(t;\var,\psi))$ satisfying (2.4).
 Moreover, since only positive or non-negative solutions of (1.1) are biologically meaningful, we restrict our framework to   positive or  non-negative initial conditions.
A  vector $x\in\R^n$ is said to be {\it positive}, or {\it non-negative}, if all its components are positive, or non-negative, respectively, and we write $x>0$, $x\ge 0$, respectively.  We define and denote  in a similar way {\it positive} and {\it non-negative functions} in $BC$, and  {\it positive} and {\it non-negative matrices} as well. As usual, we use the notation  $\R^n_+=\{x\in\R^n:x\ge 0\}$. In the space $UC_g$, a vector $c$  is identified with the constant function $\psi(s)=c$ for $s\le 0$.

 Consider the positive cone $BC^+=BC^+(\R^{2n})=\{ (\var,\psi)\in BC: \var(s),\psi(s)\ge 0$ for all $s\le 0\}$. As  set of admissible initial conditions for (1.1), we take the  subset $BC_0^+$ of $BC^+$, $BC^+_0=\{ (\var,\psi)\in BC^+:\var(0)>0,\psi(0)>0\}$. It is easy to see that all the coordinates of solutions with initial conditions in $BC^+$, respectively  $BC_0^+$, remain nonnegative, respectively positive, for all $t\ge 0$ whenever they are defined.

In the sequel, we shall consider norms $|\cdot|_d$  in $\R^N$ ($N=n$ or $N=2n$) given by $|(x_1,\dots,x_N)|_d=\displaystyle \max_{1\le i\le N}d_i|x_i|$, for some $d=(d_1,\dots,d_N)>0$. For such norms  in $\R^N$, in order to be more explicit, we denote the norm in $UC_g$  by $\|\cdot\|_{g,d}$,
$$\|\psi\|_{g,d}=\sup_{s\le 0}{|\psi(s)|_d\over g(s)}.
$$

\medskip

{\bf 3. Existence and global attractivity of a saturated equilibrium}


\med

In the absence of controls, the  Lotka-Volterra system reads as
$$ x_i^{\prime}(t)=x_i(t)\biggl(b_i-\mu_ix_i(t)-\sum_{j=1}^n a_{ij} \int_0^{\infty}K_{ij}(s)x_j(t-s)\, ds\biggr),
\eq(3.1)$$
for which
$$M_0=N+A,  \q {\rm  where}\q N=diag\ (\mu_1,\dots,\mu_n),\q A=[a_{ij}],
\eq(3.2)$$
is designated as the  interaction {\bf community matrix}. As for ordinary differential equations (ODEs), the algebraic properties of $M_0$ determine many features of  the asymptotic behaviour of solutions to (3.1) (cf.~e.g.~[4,~5,~14]). Clearly, the introduction of controls might change  the dynamics of (3.1). Here, the main aim is to use the controls to change the position of a globally attractive equilibrium, and give general criteria for its attractivity. 

For (1.1), we  define the {\bf controlled community matrix} as
$$M=N+A+C,\q {\rm  where}\q  C =diag\ \Big({{c_1d_1}\over {e_1}},\dots,{{c_nd_n}\over {e_n}}\Big).
\eq(3.3)$$
We also consider the matrices
$$\hat M_0=N-|A|,\q \hat M=N-|A|-C,\q {\rm where}\q |A|=\Big[ |a_{ij}|\Big ].
\eq(3.4)$$

Note that $(x^*,u^*)\in \R^n\times\R^n$ is an equilibrium of (1.1) if and only if
$$x_i^*=0\ {\rm or}\ (M x^*)_i=b_i, \q {\rm and} \q u_i^*={{d_i}\over {e_i}}x_i^*,\q i=1,\dots,n.$$

Throughout the paper, we shall use the definition of a {\it saturated equilibrium}.

\med

\med

 {\bf Definition 3.1}.  Let $(x^*,u^*)=(x_1^*,\dots,x_n^*,u_1^*,\dots, u_n^*)$ be an equilibrium of (1.1). We say that $(x^*,u^*)$ is a {\bf saturated equilibrium} if $(x^*,u^*)$ is non-negative and
 $$(M x^*)_i\ge b_i\q {\rm whenever}\q x_i^*=0,\q i=1,\dots,n.$$
 
 \med
 
 \noindent {\bf Remark 3.1}. We observe that if $(x^*,u^*)\ge 0$ is an equilibrium of (1.1) on the border of the positive cone $\R^n_+\times\R^n_+$,  i.e., $x_i^*=u_i^*=0$  for some $i$, and $(x^*,u^*)$ is not saturated, then $(x^*,u^*)$ is unstable. In fact,
system (1.1) and the ODE system in $\R^n_+\times \R^n_+$
$$\system{
\displaystyle x_i^{\prime}(t)&=x_i(t)\Big(b_i-\mu_ix_i(t)-\sum_{j=1}^n a_{ij} x_j(t)-c_i u_i(t)\Big), \cr
\displaystyle u_i^{\prime}(t)&=-e_i u_i(t)+d_i x_i(t), \quad i=1,2,\ldots,n , \cr}
\eq(3.5)$$
share the same equilibria. Since $\R^n_+\times \R^n_+$ is forward invariant for (3.5),
if $(x^*,u^*)\ge 0$ is an equilibrium of (3.5)  and $(x^*,u^*)$ is not saturated, then $(x^*,u^*)$ is unstable, since the characteristic equation for the linearised equation about $(x^*,u^*)$ has an eigenvalue with positive real part (cf.~e.g.~[14]).
 
 \med
 
 When analysing (1.1), our concepts of  attractivity and stability always refer to the set of {\it admissible solutions}, i.e., to solutions  $(x(t), u(t))=(x(t;\var,\psi), u(t;\var,\psi))$ with $(\var,\psi)$ in the set of admissible initial conditions. In particular,  an equilibrium $(x^*,u^*)$ of (1.1) is  {\it globally attractive} if  all  solutions $(x(t),u(t))$ of (1.1)  with initial conditions $(x_0,u_0)=(\var,\psi)\in BC^+_0, $ satisfy  $\lim_{t\to \infty} x(t)=x^*, \lim_{t\to \infty} u(t)=u^*$; and it is
 {\it globally asymptotically stable} (GAS) if it is  stable and globally attractive.
 
 We recall some concepts from matrix theory which will be used in the next sections.

 \med

 {\bf Definition 3.2}. Let $B=[b_{ij}]$ be an ${n\times n}$ matrix. 
 We say that 
   $B$ is  an {\bf  M-matrix} (respectively  {\bf non-singular M-matrix}) if $b_{ij}\le 0$ for $i\ne j$ and all its eigenvalues have non-negative (respectively positive) real parts. 
 The matrix $B$ is said to be  a {\bf P-matrix} if all its principal minors are positive.
 
 \med
   
\noindent {\bf Remark 3.2}. It is well-known that there are several equivalent ways of defining M-matrices, non-singular M-matrices and P-matrices; in [7], these matrices are also designated by {\it matrices of classes} $K_0$, $K$ and $P$, respectively. See  [1,~7,~14] for further properties of these matrices. In particular, we recall that a square matrix  
with non-positive off-diagonal entries is an M-matrix (respectively, a non-singular M-matrix) if and only if all its principal minors are non-negative (respectively, positive); so any non-singular M-matrix is a P-matrix. A related  concept is the notation of a  {\it Volterra-Lyapunov stable} ({\it VL-stable} for short) matrix, i.e.,  an ${n\times n}$ matrix $B=[b_{ij}]$ for which 
     there exists a positive vector $d=(d_1,\dots,d_n)$ such that
   $\sum_{i,j=1}^n x_id_ib_{ij}x_j<0\ {\rm for\ all}\  x=(x_1,\dots,x_n)\ne 0.$
 If $-B$ is VL-stable then $B$ is also a P-matrix;   the converse is  true for the particular case of a $2\times 2$ matrix, but not for higher dimensions. For Lotka-Volterra ODE systems of the form $x'_i= x_i[b_i-\sum_{j=1}^n a_{ij}x_j],\, 1\le i\le n$, it is known that if $-[a_{ij}]$ is VL-stable, then there is one globally stable saturated equilibrium [14, p.~199].

 \med

Consider both the original and the controlled community matrices $M_0,M$, as well as the matrices $\hat M_0,\hat M$, cf.~(3.2)--(3.4).  For the uncontrolled system (3.1), it was shown in [4, Corollary 4.1] that  if $\hat M_0$ is a non-singular M-matrix, then there is a unique saturated equilibrium  of (3.1),
 which is a global attractor of all solutions with initial conditions $x_0=\var \in BC_0^+(\R^n)$.
  The idea now is to prove a similar result for system (1.1). We start by studying  the existence of a saturated equilibrium and the boundedness of solutions  to (1.1).

\proclaim{Theorem 3.1}.  Assume that $M$ is a P-matrix, where $M$ is  the controlled community matrix    in (3.3). Then,  there is a unique saturated equilibrium $(x^*,u^*)$ of (1.1).

 {\it Proof}. If $M$ is a P-matrix, then  for each vector $b=(b_1,\dots,b_n)\in\R^n$ there is a unique non-negative vector $x^*$ such that $M x^*\ge b$ and $(M x^*)_i= b_i$ if $x_i^*>0$ [1, p.~274]. With $u^*=(u_1^*,\dots, u_n^*)$ where $u_i^*={{d_i}\over {e_i}}x_i^*$, this means that $(x^*,u^*)$ is the unique saturated equilibrium of (1.1).\ter

If all coefficients in (1.1) are positive, then clearly all positive solutions are bounded, since the inequalities $x_i'(t)\le x_i(t) (b_i-\mu_ix_i(t))$ hold, and positive solutions of the logistic ODEs $y'(t)=y(t) (b_i-\mu_iy(t))$ are bounded. This is not however the case if we allow some of the coefficients $a_{ij}$ to be  negative, unless further constraints on $M_0$ are imposed.

\proclaim{Lemma 3.1}. Assume that the matrix $\hat M_0$ in (3.4) is a non-singular M-matrix. Then, all  solutions of (1.1) with initial conditions (1.2) are defined and bounded on $[0,\infty)$.

{\it Proof}. Solutions of (1.1) with initial conditions (1.2) are positive, whenever they are defined. For (1.1) written in the abstract form $X'(t)=F(X_t)$, the function $F$ transforms bounded sets of $UC_g(\R^{2n})$ into bounded sets of $\R^{2n}$, hence solutions are defined  on compact intervals $[0,\al], \forall \al>0$, and therefore   on $[0,\infty)$. 

Since $\hat M_0$ is a non-singular M-matrix, there is a positive vector $\eta=(\eta_1,\dots,\eta_n)$ such that $\hat M_0\eta>0$ [7], i.e.,
 $$\mu_i\eta_i>\sum_{j=1}^n |a_{ij}|\eta_j,\q i=1,\dots,n.$$
Choose an arbitrarily small $\de>0$  so that
$$\mu_i-(1+\de)\sum_{j=1}^n |a_{ij}| {{\eta_j}\over {\eta_i}}>0,\q i=1,\dots,n,
 \eq(3.6)$$
 and  a function $g$ for which (g1)-(g3) and (2.2) hold. 
 
 For $\bar\eta =(\eta_1^{-1},\dots ,\eta_n^{-1},e_1(d_1\eta_1)^{-1},\dots, e_n(d_n\eta_n)^{-1})$, we further
 consider  $\R^{2n}$ equipped with the norm $|\cdot|_{\bar\eta}$  given by
$$|(x_1,\dots,x_n,u_1,\dots,u_n)|_{\bar\eta}=\max \left \{ \max_i({1\over \eta_i }|x_i|),\max_i({{e_i}\over {\eta_i d_i}}|u_i|)\right\}.$$

Let  $(x(t),u(t))=(x_1(t),\dots,x_n(t),u_1(t),\dots,u_n(t))$ be a  positive solution of (1.1). We claim that
$$\sup_{t\ge 0} |(x(t),u(t)) |_{\bar \eta}<\infty.
\eq(3.7)$$

For the sake of contradiction, assume that (3.7) fails. Then, for any $K>0$ there exists $T>0$ such that 
 $$|(x(T),u(T))|_{\bar\eta}\ge | (K,\dots,K)|_{\bar \eta}\q {\rm and}\q |(x(T),u(T))|_{\bar\eta}\ge |(x(t),u(t))|_{\bar\eta},\ 0\le t\le T.\eq(3.8)$$
  Consider (3.8) with $K$ such  that $| (K,\dots,K)|_{\bar \eta}> \| (x_0,u_0)\|_{g,\bar \eta}=
\displaystyle \sup_{s\le 0}{{  |(x(s),u(s)) |_{\bar \eta}}\over {g(s)}}$, and sufficiently large to be specified later.
 
 \med
 
 If $|(x(T),u(T))|_{\bar\eta}={{ e_i}\over {\eta_id_i}}u_i(T)>{1\over {\eta_i}} x_i(T)$ for some $i\in\{1,\dots,n\}$, from
  (1.1) we obtain
 $$u_i'(T)< -e_i u_i(T)+d_i {{e_i}\over {d_i}} u_i(T)\le 0,$$
 which is not possible since the definition of $T$ implies $u_i'(T)\ge 0$.
  Thus, $|(x(T),u(T))|_{\bar\eta}={1\over {\eta_i}} x_i(T)$
 for some $i\in\{1,\dots,n\}$. Clearly $x_i'(T)\ge 0$. 
 
 Let $0\le t\le T$ and $s\le 0$. Note that ${1\over {\eta_j}}{{x_j(t-s)}\over {g(-s)}}\le {1\over {\eta_j}}{{x_j(t-s)}\over {g(t-s)}}< {1\over {\eta_j}}K\le{1\over {\eta_i}}x_i(T)$ if $t-s\le 0$, and ${1\over {\eta_j}}{{x_j(t-s)}\over {g(-s)}}\le{1\over {\eta_i}}x_i(T)$ if $0\le t-s\le T$. 
 Hence
$$\eqalign{
\eta_j^{-1}\left | a_{ij}\int_0^{\infty}K_{ij}(s)x_j(t-s)\, ds\right| &\le \eta_j^{-1}|a_{ij}|\int_0^{\infty}g(-s)K_{ij}(s){{x_j(t-s)}\over {g(-s)}}\, ds\cr
&< (1+\de) \eta_i^{-1}|a_{ij}|x_i(T),\ j=1,\dots,n.\cr}
\eq(3.9)$$
From (1.1) and (3.9), we obtain
$$
0\le x_i'(T) \le x_i(T)\bigg[b_i-\Big (\mu_i-(1+\de)\sum_{j=1}^n |a_{ij}| {{\eta_j}\over {\eta_i}}\Big )x_i(T)\bigg].$$
By (3.6), this is a contradiction   if $K$ is chosen so that
$K>b_i \Big (\mu_i-(1+\de)\sum_{j=1}^n |a_{ij}| {{\eta_j}\over {\eta_i}}\Big )^{-1}.$\ter

In fact,  a better criterion for the uniform boundedness of all positive solutions of (1.1) will be given later (cf. Theorem 3.3).

 Note that  $\hat M_0=\hat M+ C$ where   $C$ is a positive diagonal matrix. By Theorem 5.1.1 of  [7], it follows that if  $\hat M$  an M-matrix, then
  $\hat M_0$ is a non-singular M-matrix. Now, if $\hat M_0$ is a non-singular M-matrix, 
there is a positive vector $\eta=(\eta_1,\dots,\eta_n)$
such that $\hat M_0\eta>0$, i.e., $\mu_i\eta_i>\sum_{j=1}^n |a_{ij}|\eta_j$ for $1\le i\le n$ 
 (cf.~[7]); in particular, this implies 'diagonal dominance' of $M_0$, in the sense that
 $(\mu_i+a_{ii})\eta_i>\sum_{j\ne i} |a_{ij}|\eta_j$ for $1\le i\le n$. From [14, p.~201],  it follows that  if $\hat M_0$  is a non-singular M-matrix, then $-M_0$ (and hence $-M$ as well) is
 VL-stable, and therefore a P-matrix.

  \proclaim{Lemma 3.2}. Assume that the matrix $\hat M$ in (3.4) is an M-matrix, and that the unique saturated equilibrium $(x^*,u^*)$ of (1.1) is positive. Then $(x^*,u^*)$ is locally asymptotically stable.
  
  {\it Proof}. As observed, if $\hat M$ is an M-matrix then $M$ is a P-matrix, and there is a unique saturated equilibrium $(x^*,u^*)$. The linearisation of (1.1) about $(x^*,u^*)$ is  given by
  $$\left[\matrix {y'(t)\cr v'(t)\cr}\right]=-\left ( B\left[\matrix {y(t)\cr v(t)\cr}\right]+{\cal L}\left[\matrix {y_t\cr v_t\cr}\right]\right),
  \eq(3.10)$$
 with $y(t),v(t)\in \R^n$, and the $(2n)\times (2n)$ matrix $B$ and the linear operator ${\cal L}:BC_g(\R^{2n})\subset UC_g(\R^{2n})\to \R^{2n}$ defined below:
 $$B=diag (\al_1,\dots,\al_n,e_1,\dots,e_n),\q {\cal L}=({\cal L}_1,\dots,{\cal L}_{2n}),$$
 where
 $$
 \al_i=\system{&\mu_ix_i^*,& \ {\rm if}\ x_i^*>0,\cr &\sum_{j=1}^na_{ij}x_j^*-b_i,& \ {\rm if}\ x_i^*=0,\cr},$$ and 
$$
\system{&{\cal L}_{i}(\var,\psi)=x_i^*\sum_{j=1}^n a_{ij} \int_0^{\infty}K_{ij}(s)\var_j(-s)\, ds+x_i^*
 c_i \int_0^{\infty}G_i(s)\psi_i(-s)\, ds,\cr
& {\cal L}_{n+i}(\var,\psi)=-d_i\var_i(0),\qq i=1,\dots,n,\cr}$$
 for $(\var,\psi)=(\var_1,\dots,\var_n,\psi_1,\dots,\psi_n)$. Note that $\al_i\ge 0$ for $1\le i\le n$.
 
 For $(\vec e_1,\dots,\vec e_{2n})$ the canonical basis of $\R^{2n}$, define
 $L:=B+\Big [{\cal L}_i(\vec e_j)\Big ]_{i,j}$ and $ \hat L:=B-\Big [|{\cal L}_i(\vec e_j)|\Big ]_{i,j}\, (1\le i,j\le 2n)$. We have
 $$\Big [{\cal L}_i(\vec e_j)\Big ]_{i,j=1}^n=\left [\matrix {A(x^*)&C(x^*)\cr -D&0\cr}\right ],$$
 with $A(x^*)=\big [x_i^*a_{ij}\big ]_{i,j}\, (1\le i,j\le n), C(x^*)=diag\, (x_1^*c_1,\dots,x_n^*c_n)$ and $D=diag\, (d_1,\dots,d_n)$. Now suppose that $x^*>0$. It is easy to see that the matrices $L,\hat L$ are equivalent to, respectively,
  $$\left [\matrix {M(x^*)&C(x^*)\cr 0&E\cr}\right ],
  \q \left [\matrix {\hat M(x^*)&-C(x^*)\cr 0&E\cr}\right ],$$ 
 where $E=diag\, (e_1,\dots,e_n), N(x^*)=diag\, (\mu_1x_1^*,\dots,\mu_nx_n^*), \tilde C(x^*)=diag\, (x_1^*{{c_1d_1}\over {e_1}},\dots,x_n^*{{c_nd_n}\over {e_n}}),$
 $M(x^*)=N(x^*)+A(x^*)+\tilde C(x^*)$ and $ \hat M(x^*)=N(x^*)-A(x^*)-\tilde C(x^*).$ Note that $M(x^*),\hat M(x^*)$ are obtained from $M,\hat M$, respectively, by multiplying each line $i$ by $x_i^*$.
 Hence, it follows that $\det L\ne 0$ and that $\hat L$ is an M-matrix as well. From  [4], we derive that the linear system (3.10) is exponentially asymptotically stable.\ter
 
 We remark that if $x_i^*=0$, then the $i$th-line of the above matrix $\Big [{\cal L}_i(\vec e_j)\Big ]$ is zero. Hence, a  saturated equilibrium $(x^*,u^*)$  of (1.1) on the boundary of the positive cone is not necessarily asymptotically stable, thus although its linearisation (3.10) is stable, one cannot deduce that $(x^*,u^*)$ is stable as a solution of (1.1).
 
Our main general result on the global attractivity of the saturated equilibrium is given below.
  
  \proclaim{Theorem 3.2}. Assume that the matrix $\hat M$ in (3.4) is an M-matrix. Then, there is a unique saturated equilibrium $(x^*,u^*)$ of (1.1),
 which is a global attractor of all solutions with initial conditions (1.2). Moreover, if $x^*>0$, then $(x^*,u^*)$ is GAS.

{\it Proof}. Since $\hat M$ is an M-matrix,  from Theorem 3.1 and Lemma 3.1 we  conclude that there is a unique saturated equilibrium $(x^*,u^*)$ of (1.1) and that all positive solutions are defined and bounded on $[0,\infty)$. Lemma 3.2 shows that  $(x^*,u^*)$ is stable if it is a positive equilibrium.
We now need to show that $(x^*,u^*)$ is a global attractor of all positive solutions of (1.1). 

\smal

Denote $I_n$  the $n\times n$ identity matrix. If $\hat M$ is an M-matrix, then for any $\de_0>0$ the matrix $\de_0 I_n+\hat M$  is a non-singular M-matrix.  Fix any $\de_0>0$
and a positive vector $\eta=(\eta_1,\dots,\eta_n)$
such that $(\de_0 I_n+\hat M)\eta>0$, i.e.,$$(\de_0 +\mu_i-c_i{{d_i}\over {e_i}})\eta_i>\sum_{j=1}^n |a_{ij}|\eta_j,\q i=1,\dots,n.
$$ 
Choose $\de>0$ such that  
$$(\de_0 +\mu_i-c_i{{d_i}\over {e_i}})\eta_i-(1+\de)\sum_{j=1}^n |a_{ij}| \eta_j>0,\q 1\le i\le n,
\eq(3.11)$$
 and  a function $g$ for which conditions (g1)-(g3) and (2.2) are fulfilled. We abuse the notation, and denote both  norms  in $\R^{2n}$ and  in $\R^{n}$ by $|\cdot|_{\bar\eta}$, where
$$|(x_1,\dots,x_n,u_1,\dots,u_n)|_{\bar\eta}:=\max _{1\le i\le n} \left \{ \max \Big ({1\over \eta_i }|x_i|,{{e_i}\over {\eta_i d_i}}|u_i|\Big )\right\}\q {\rm in}\q \R^{2n},$$
 $$|x|_{\bar\eta}:=\max _{1\le i\le n}{1\over \eta_i }|x_i|\q {\rm in}\q \R^n,$$
  and consider $UC_g(\R^{2n})$, $UC_g(\R^n)$ equipped with the norms $\|\cdot\|_{g,\bar \eta}$.
 
\med

Let $(x(t),u(t))$ be a positive solution of (1.1).  With the change of variables
$$y_i(t)=x_i(t)-x_i^*,\ v_i(t)=u_i(t)-u_i^*,\q i=1,\dots,n,$$
system (1.1) together with definition 3.1 lead to
$$\eqalignno{
 y_i^{\prime}(t)&=-(y_i(t)+x_i^*)\biggl(\mu_iy_i(t)+\sum_{j=1}^n a_{ij} \int_0^{\infty}K_{ij}(s)y_j(t-s)\, ds\cr
 &\hskip 3cm+c_i \int_0^{\infty}G_i(s)v_i(t-s)\, ds\biggr),\ {\rm if}\ x_i^*>0,&(3.12) \cr
y_i^{\prime}(t)&\le-y_i(t)\biggl(\mu_iy_i(t)+\sum_{j=1}^n a_{ij} \int_0^{\infty}K_{ij}(s)y_j(t-s)\, ds\cr
 &\hskip 3cm+c_i \int_0^{\infty}G_i(s)v_i(t-s)\, ds\biggr),\ {\rm if}\ x_i^*=0,&(3.13) \cr
 v_i^{\prime}(t)&=-e_i v_i(t)+d_i y_i(t), \quad i=1,2,\ldots,n .&(3.14)\cr}
$$
 Define
$$\liminf_{t\to\infty} y_i(t)=-l_i,\q \limsup_{t\to\infty} y_i(t)=L_i,\q i=1,\dots,n,
\eq(3.15)$$
and set
$$l=\max_{1\le i\le n}{{l_i}\over {\eta_i}},\q L=\max_{1\le i\le n}{{L_i}\over {\eta_i}},\q U=\max(l,L).$$
 Integrating (3.14), we get
$$v_i(t)=v_i(0)e^{-e_it}+{d_i}e^{-e_it}\int_0^t e^{e_is}y_i(s)\, ds,\q t\ge 0,
\eq(3.16)$$
and therefore
 $$-x_i^*\le -l_i\le {{e_i}\over {d_i}} \liminf_{t\to \infty}v_i(t)
 \le {{e_i}\over {d_i}} \limsup_{t\to \infty}v_i(t)\le \L_i< \infty.
 \eq(3.17)$$
 
Since $U\ge 0$, it is enough to prove that $U=0$. In order to get a contradiction, assume $U>0$.

Denote $I=\{1,\dots,n\}, I_1=\{ i\in I: \eta_i^{-1}L_i=U\}$ and $ I_2=\{ i\in I: \eta_i^{-1}l_i=U\}$. The assumption $U>0$ implies that $x_i^*>0$ if $i\in I_2$; otherwise, with $x_i^*=0$ we get $\liminf_{t\to\infty} x_i(t)=\liminf_{t\to\infty} y_i(t)=-l_i=-\eta_iU\ge 0$, thus $U=0$.

The coordinates $ y_j(t), v_j(t)$ are uniformly bounded for $t\ge 0$, thus,
as remarked in Section 2, the positive orbit $\{ (y_t,v_t):t\ge 0\}$ is precompact in $UC_g(\R^{2n})$. 

Take  any  sequence $(t_k)$ with $t_k\to\infty$. Thus, there is a subsequence of $(y_{t_k}, v_{t_k})$, still denoted by  $(y_{t_k}, v_{t_k})$, converging  to some $(\phi,\psi)$ in $UC_g(\R^{2n})$. Let $\phi_j,\psi_j\, (1\le j\le n)$ be the components of $\phi,\psi$, respectively.
Take any $\vare>0$ and let $t^*>0$ be such that $\eta_j^{-1}|y_j(t)|\le U+\vare,$ for $t\ge t^*,1\le j\le n$.   For any $s\ge 0$, if $k$ is large so that $t_k-s\ge t^*$, then
$\eta_j^{-1}|y_{t_k,j}(-s)|/g(-s)|=\eta_j^{-1}|y_j(t_k-s)|/g(-s)\le \eta_j^{-1} |y_j(t_k-s)|\le U+\vare$, therefore we get
$\eta_j^{-1}\|\phi_j\|_g\le U+\vare$.
In a similar way, 
we obtain 
$\eta_j^{-1}{{e_i}\over {d_i}}\|\psi_{j}\|_g\le U+\vare$. Hence, we conclude that 
$\|(\phi,\psi)\|_{g,\bar \eta}\le U$. Moreover, if  $i\in I_1\cup I_2$ and $(t_k)$ is chosen in such a way that $\eta_i^{-1}|y_i(t_k)|\to U$, we further deduce that $\|\phi\|_{g,\bar \eta}=\eta_i^{-1}|\phi_i(0)|= U$ and that $y_{t_k,j},v_{t_k,j}$  converge uniformly  to $\phi_j,\psi_j$, respectively, on each compact set  of $[0,\infty)$.

Fix $i\in I_1\cup I_2$. By the fluctuation lemma, take a sequence $(t_k)$ with $t_k\to \infty, y'_i(t_k)\to 0$ and $\eta_i^{-1}y_i(t_k)\to \system{U, &\ {\rm if}\ i\in I_1\cr -U, &\ {\rm if}\ i\in I_2\cr}$. As above, we may assume that $(y_{t_k},v_{t_k})\to (\phi,\psi)\in UC_g(\R^{2n})$ for the norm $\|\cdot\|_{g,\bar \eta}$.

First, we consider the case $i\in I_1$, thus  $\eta_i^{-1} y_i(t_k)\to U$.

Since the linear operator $\psi\mapsto \int_0^\infty  G_i(s)\psi(-s)\, ds$, defined for $ \psi\in BC(\R)\subset UC_g(\R)$,  is bounded, there exists
$$\nu:=\lim_{k\to\infty} \int_0^\infty G_i(s)v_i(t_k-s)\, ds=\int_0^\infty G_i(s)\psi_i(-s)\, ds.
$$
From (3.17), we have $|\psi_i (-s)|\le  \eta_i{{d_i}\over {e_i}}U$ for any $s\ge 0$, thus $\nu\ge - \eta_i{{d_i}\over {e_i}}U$. From hypothesis  (1.3),  the inequality
 $\nu= - \eta_i{{d_i}\over {e_i}}U$ implies that $\psi_i(0)= - \eta_i{{d_i}\over {e_i}}U$. But,
from (3.16) we have
$$v_i(t_k)=v_i(0)e^{-e_it_k}+{d_i}\int_0^{t_k} e^{-e_iu}y_{t_k,i}(-u)\, du,$$
and from the Lebesgue's dominated convergence theorem it follows that
 $$\lim_k v_i(t_k)=\psi_i(0)=d_i\int_0^{\infty} e^{-e_iu} \phi_i(-u)\, du.$$
Since $\phi_i$ is a continuous function with $\eta_i^{-1}|\phi_i(-s)|\le U$ for  $s>0$ and $\eta_i^{-1}\phi_i(0)=U$, then
$\eta_i^{-1}\int_0^{\infty} e^{-e_is} \phi_i(-s)\, ds>-U/e_i.$
We therefore conclude that
$$\eta_i^{-1}\nu>-{{d_i}\over {e_i}}U.\eq(3.18)$$
Moreover, in spite of the use of a specific vector $\eta=\eta(\de_0)$ and norm  $\|\cdot\|_{g,\bar \eta}$ in $UC_g(\R^n)$, obviously the limit $\nu$ does not depend on the chosen norm $|\cdot|_{\bar\eta}$ in $\R^n$.

Next, denote 
$$H_i(t)=\mu_iy_i(t)+\sum_{j=1}^n a_{ij} \int_0^{\infty}K_{ij}(s)y_j(t-s)\,ds+c_i \int_0^{\infty}G_i(s)v_i(t-s)\, ds.
\eq(3.19)$$
 From (3.12)-(3.13), we obtain
 $$y'_i(t_k)\le -(y_i(t_k)+x_i^*)H_i(t_k).$$
From (2.2), we have (cf. (3.9))
$$
\left | a_{ij}\int_0^{\infty}K_{ij}(s)y_j(t_k-s)\, ds\right|\le
|a_{ij}| \int_0^{\infty}g(-s)K_{ij}(s){{|y_j(t_k-s)|}\over {g(-s)}}\, ds
\le (1+\de)|a_{ij}|\|y_{t_k,j}\|_g
\eq(3.20)$$
and  this  leads to
$$\eqalign{
H_i(t_k)&\ge \mu_iy_i(t_k)- (1+\de)\sum_{j=1}^n |a_{ij}|\,  \|y_{t_k,j}\|_{g,\bar\eta} +c_i \int_0^\infty G_i(s)v_i(t_k-s)\, ds\cr
&\ge  \mu_iy_i(t_k)- (1+\de)\sum_{j=1}^n |a_{ij}|\, \eta_j \|y_{t_k}\|_{g,\bar\eta} +c_i \int_0^\infty G_i(s)v_i(t_k-s)\, ds.\cr}\eq(3.21)
$$
By letting $k\to \infty$, from (3.11) and (3.21) we have
$$0\ge  \bigg (\mu_i\eta_i-(1+\de)\sum_{j=1}^n |a_{ij}|\eta_j\bigg)U+c_i\nu\ge \Big (c_i{{d_i}\over {e_i}}-\de_0\Big )\eta_iU+c_i\nu .\eq(3.22)
$$
Since $\de_0>0$ is arbitrarily small, this yields $\nu\le-{{d_i}\over {e_i}}\eta_iU$,  which is not possible in view of (3.18).

\med 

Now, consider the case  $i\in I_2$. Then,  $\eta_i^{-1}y_i(t_k)\to -U$ and (3.12) holds.

If $y_i(t)$ is eventually monotone, then $y_i(t)\to -\eta_iU$ and $v_i(t)\to -{{d_i}\over {e_i}}\eta_iU$. 
Using arguments similar to the ones above,  we obtain 
$$\eqalign{
H_i(t_k)\le \mu_iy_i(t_k)&+(1+\de)\sum_{j=1}^n |a_{ij}|\,  \eta_j \|y_{t_k}\|_{g,\bar\eta}+c_i \int_0^\infty G_i(s)v_i(t_k-s)\, ds\cr
&\to \Big[-(\mu_i+c_i{{d_i}\over {e_i}})\eta_i+(1+\de)\sum_{j=1}^n |a_{ij}| \eta_j \Big ]U<0.}
\eq(3.23)$$
Since $y'_i(t_k)=-(y_i(t_k)+x_i^*)H_i(t_k)$, using the above estimate  we obtain
$$0\ge (-\eta_iU+x_i^*)\bigg[(\mu_i+c_i{{d_i}\over {e_i}})\eta_i-(1+\de)\sum_{j=1}^n |a_{ij}| \eta_j\bigg]U,$$
 and thus
 $y_i(t)\to -x_i^*=-\eta_iU$ as $t\to\infty$.  Since $y_i(t)>-x_i^*$ for $t>0$, this is only possible if $y'_i(t)\le 0$ for $t$ large, so that $\eta_i^{-1}y_i(t)\searrow -U$.
But in this case from (3.12) it follows   that $ H_i(t)\ge 0$ for $t$ large, which contradicts $(3.23)$.

If  $y_i(t)$ is not eventually monotone, then we can assume that $ y_i(t_k)$ is a sequence of minima, so that $H_i(t_k)=0$, and
 this case is treated as the case $i\in I_1$. These arguments show that $U=0$, and the proof is complete.\ter 
 \med 
 
 \noindent {\bf Remark 3.3}.  As referred to in the introduction,  clearly  the above proof applies to systems (1.4). In fact, with the terms $a_{ij} \int_0^{\infty}K_{ij}(s)x_j(t-s)\, ds, c_i \int_0^{\infty}G_i(s)u_i(t-s)\, ds$  replaced by  the more general linearities 
$$a_{ij} \int_0^{\infty}x_j(t-s)\, d\eta_{ij}(s),\q c_i \int_0^{\infty}u_i(t-s)\, d\nu_i(s),$$
where $\eta_{ij},\nu_i$ are normalized bounded variation functions and $\nu_i$ are non-decreasing, we use (2.3) instead of (2.2), 
 the estimates (3.20) are replaced by
  $$\eqalign{
\left | a_{ij}\int_0^{\infty}y_j(t_k-s)\, d\eta_{ij}(s)\right| &\le |a_{ij}|\int_0^{\infty}g(-s){{|y_j(t_k-s)|}\over {g(-s)}}\, d|\eta_{ij}(s)|\cr
&\le (1+\de) |a_{ij}|\|y_{t_k,j}\|_g,\q i,j=1,\dots,n,\cr}
$$
the limit $\nu$ is now given by $\nu=\int_0^\infty \psi_i(-s)\, d\nu_i(s)$, and all the other arguments are valid.

\med

With the usual notation of
$$a_{ij}=a_{ij}^+-a_{ij}^-,\q
{\rm where}\q
a_{ij}^+=\max \{a_{ij},0\},\ 
a_{ij}^-=\max\{-a_{ij},0\},$$
we denote 
$$M_0^-=diag\ (\mu_1,\dots,\mu_n)-A^-\, ,\q {\rm where}\q A^-=\big [a_{ij}^-\big ].
 \eq(3.24)$$
Note that $M_0^-\ge \hat M_0$, hence in general imposing that $M_0^-$ is a non-singular M-matrix is weaker than requiring that $\hat M_0$ is a non-singular M-matrix. We  now give sufficient conditions for the dissipativeness of (1.1), improving Lemma 3.1.

 \proclaim{Theorem 3.3}.  If $M_0^-$ is a non-singular M-matrix, then
  (1.1) is dissipative; i.e., there exists $K>0$ such that $\limsup_{t\to\infty} x_i(t)\le K, \limsup_{t\to\infty} u_i(t)\le K,\, 1\le i\le n$, for all solutions $(x(t),u(t))$ of (1.1) with initial conditions (1.2).

{\it Proof}. A   solution $(x(t),u(t))$ of (1.1) with  initial condition  $(x_0,u_0)=(\var,\psi)\in BC^+_0$ satisfies 
$$\system{
\displaystyle x_i^{\prime}(t)&\le x_i(t)\biggl(b_i-\mu_ix_i(t)+\sum_{j=1}^n a_{ij}^- \int_0^{\infty}K_{ij}(s)x_j(t-s)\, ds\biggr), \cr
\displaystyle u_i^{\prime}(t)&=-e_i u_i(t)+d_i x_i(t), \quad i=1,2,\ldots,n . \cr}
$$
Let  $(X(t),U(t))$ be  the solution of the  system
$$\system{
\displaystyle X_i^{\prime}(t)&=X_i(t)\biggl(b_i-\mu_iX_i(t)+\sum_{j=1}^n a_{ij}^- \int_0^{\infty}K_{ij}(s)X_j(t-s)\, ds\biggr) \cr
\displaystyle U_i^{\prime}(t)&=-e_i U_i(t)+d_i X_i(t), \quad i=1,2,\ldots,n . \cr}
\eq(3.25)$$
with the initial conditions $X_0=\var, U(0)=\psi (0)$. Since (3.25) is cooperative, or in other worths, it satisfies the quasi-monotonicity condition  in [22, Chapter 5], by comparison results it follows that $x(t)\le X(t),u(t)\le U(t)$. From [4, Corollary 4.1],  $(X(t),U(t))\to (X^*,U^*)$ as $t\to\infty$, where $(X^*,U^*)$ is the saturated equilibrium of (3.25). Thus, the solutions $(x(t),u(t))$ of the initial value problems (1.1)-(1.2) satisfy $\limsup_{t\to\infty}x_i(t)\le X_i^*,
 \limsup_{t\to\infty}u_i(t)\le U_i^*,\, 1\le i\le n$.\ter
\smal

 Our setting contemplates all the possibilities for the signs of the coefficients $b_i,a_{ij}$ in (3.1). 
 In biological terms, the most interesting cases are however: (i) $a_{ij}\ge 0$ for $i\ne j$ (competitive systems); 
 (ii) $a_{ij}\le 0$ for $i\ne j$ (cooperative systems); (iii) $a_{ij}>0,a_{ji}<0$ (predator-prey systems) if species  $i$ is a prey for the predator species $j$, $i\ne j$.
 On the other hand, the existence of a positive equilibrium depends heavily on the coefficients $b_i$'s, and  can be studied in more detail by using Cramer's rule.  Nevertheless,   a criterion for cooperative systems is given here.
 
 \proclaim{Theorem 3.4}. Consider (1.1) with $b_i>0$ and $a_{ij}\le 0$ for all $i\ne j$. If $M$ is a non-singular M-matrix, then there exists a unique positive equilibrium of (1.1).

 {\it Proof}.   Denote $b=(b_1,\dots ,b_n)$. Since $M$ is non-singular M-matrix, then $M^{-1}\ge 0$  [1]. This implies that $(M^{-1}b)_i=0$ if and only if the $i$th line of $M^{-1}$ is zero, which is not possible. Therefore,  $x^*:=M^{-1}b$ is a positive vector, and $(x^*,u^*)$, with $u_i^*={{d_i}\over{e_i}}x_i^*,\, 1\le i\le n$, is a positive equilibrium  of (1.1).\ter

 \bigskip

 {\bf 4.  Extinction and stability}

 \med

For the results in this section, it is important to consider systems (1.1) with $K_{ij}$ non-negative, or more general systems (1.4) with $\eta_{ij}$ non-decreasing, $i,j=1,\dots,n$. Straightforward generalisations for the situation of $K_{ij}$ in (1.1) changing signs or $\eta_{ij}$ non-monotone on $[0,\infty)$ can however be derived (cf.~[4] for the case of uncontrolled Lotka-Volterra models).

%

We now seek for better sufficient conditions for extinction of either all or part  of the populations. Together with the  controlled Lotka-Volterra system (1.1), consider the ODE system (3.5),
 and write (3.5) in the form $X'(t)=F(X(t))$ for $X(t)=(x_1(t),\dots,x_n(t),u_1(t),\dots, u_n(t))$.

Define $\la_i=\mu_i+{{c_id_i}\over {e_i}}$. If  $X^*=(x^*,u^*)=(x_1^*,\dots,x_n^*,u_1^*,\dots,u_n^*)$ is an equilibrium of (3.5), then
 $$DF(X^*)=\left[\matrix{Df(x^*)&-C(x^*)\cr D&-E}\right],$$
 where $C(x^*)=diag\, (c_1x_1^*,\dots, c_nx_n^*), D=diag\, (d_1,\dots,d_n), E=diag\, (e_1,\dots,e_n)$, and
  ${{\p f_i}\over {\p x_i}}(x^*)=b_i-(\la_i+a_{ii})x_i^*-\sum_{j=1}^n a_{ij}x_j^*-(\mu_i+a_{ii})x_i^*$,  ${{\p f_i}\over {\p x_j}}(x^*)=-a_{ij}x_i^*$ if $i\ne j$. Note that
   ${{\p f_i}\over {\p x_i}}(x^*)=-(\mu_i+a_{ii})x_i^*$ if $x_i^*>0$, otherwise ${{\p f_i}\over {\p x_i}}(x^*)\le -(\mu_i+a_{ii})x_i^*$ .

 For the trivial equilibrium, we have $C(0)=0$, hence the spectrum of $DF(0)$ is $\sigma (DF(0))=\{ b_1,\dots,b_n, -e_1,\dots,-e_n\}$. We therefore conclude that zero  is a stable equilibrium for the linearisation of  (3.5)  at zero (which is also the linearisation of  (1.1))  if and only if $b_i\le 0, 1\le i\le n$, and the introduction of the controls does not change its stability.   

Let $M$ be a P-matrix.  If $b_i\le 0, 1\le i\le n$, then 0 is the  saturated equilibrium; on reverse, if  $b_i>0$ for some  $i$, then $DF(0)$ is unstable and zero is not the saturated equilibrium. 
 In the latter case,  we have seen that (1.1) is dissipative  if $M_0^-$ is a non-singular M-matrix -- then there exists a compact global attractor ([11, Theorem 3.4.8]), which however need not  be the saturated equilibrium.
When 0 is  saturated, rather than Theorem 3.2, next result provides a better criterion for extinction of all populations.

 \proclaim{Theorem 4.1}. Assume that  $M$ is a P-matrix. The equilibrium 0 is the saturated equilibrium of (1.1) if and only if
 $b_i\le 0$ for $1\le i\le n$. In this case, if $M_0^-$ is an M-matrix, where $M_0^-$ is defined as in (3.24), then the equilibrium 0 of (1.1) is globally attractive.
 
 {\it Proof}. 
 Since $M_0^-$ is an M-matrix, for any arbitrarily small $\de_0>0$, consider a positive vector $\eta=(\eta_1,\dots,\eta_n)$ such that $(M_0^-+\de_0I_n)\eta >0$ [7].
 Let $(x(t), u(t))$ be a solution of (1.1). After a scaling $x_i\mapsto \bar x_i=\eta_i^{-1}x_i, u_i\mapsto \bar u_i=\eta_i^{-1}u_i,\ 1\le i\le n$,  and dropping the bars for the sake of simplicity, we may suppose that $(x(t), u(t))$ is a solution of (1.1) and that $(M_0^-+\de_0I_n)\eta >0$ with $\eta=(1,\dots,1)$. Next, choose $\de >0$ small and $g$ satisfying (g1)-(g3) and (2.2), with
$$\de_0+\mu_i-(1+\de)\sum_{j=1}^n a_{ij}^- >0,\q 1\le i\le n.
\eq(4.1)$$
Define $L_i=\limsup_{t\to\infty} x_i(t)$ and $U=\max_{1\le i\le n} L_i$. 
For the sake of contradiction, assume that $U>0$, and choose $i\in \{1,\dots,n\}$ such that $L_i=U$. 
Consider a sequence $(t_k)$ with $t_k\to\infty$, $x_i'(t_k)\to 0, x_i(t_k)\to U$ as $k\to\infty$. 
 We now argue as in the proof of Theorem 3.2, omitting some of the details. For some subsequence of $(x_{t_k}, u_{t_k})$, still denoted by $(x_{t_k}, u_{t_k})$, there is $(\phi,\psi)\in BC^+(\R^{2n})\subset UC_g(\R^{2n})$ such that $x_{t_k}\to \phi,  u_{t_k}\to \psi$, and $\|\phi\|_g=U=\phi_i(0)$. Next, from (1.3) and the fact that $\psi_i(0)=d_i \int_0^\infty e^{-e_iu}\phi_i(-u)\, du>0$, we obtain $\nu:=\lim_{k\to\infty}\int_0^\infty G_i(s)u_i(t_k-s)\, ds=\int_0^\infty G_i(s)\psi_i(-s)\, ds>0$. Choose $\de_0>0$ small and $k$  large so that $c_i\int_0^\infty G_i(s)u_i(t_k-s)\, ds>\de_0U$.
Since $b_i\le 0$,  for $k$  large estimates as in (3.9) yield
$$\eqalign{x_i'(t_k)&\le x_i(t_k)
\biggl(b_i-\mu_ix_i(t_k)+\sum_{j=1}^n a_{ij}^- \int_0^{\infty}K_{ij}(s)x_j(t_k-s)\, ds-c_i\int_0^\infty G_i(s)u_i(t_k-s)\, ds\biggr)\cr
\le &-x_i(t_k) \biggl(\mu_ix_i(t_k)-\sum_{j=1}^n a_{ij}^- \|x_{t_k,j}\|_g\int_0^{\infty}g(-s)K_{ij}(s)\, ds+\de_0U\biggr)
\cr
\le &-x_i(t_k) \biggl(\mu_ix_i(t_k)-(1+\de)\sum_{j=1}^n a_{ij}^- \|x_{t_k,j}\|_g+\de_0U\biggr)\, .
\cr}$$
By letting $k\to\infty$ we obtain
$$0\ge \Big [\de_0+\mu_i-(1+\de)\sum_{j=1}^n a_{ij}^-\Big ] U,$$
which contradicts (4.1). Hence $U=0$, and the proof is complete.\ter

Consider now the case of a saturated equilibrium  $(x^*,u^*)\ne 0$ of (1.1) with $x^*\in \p(\R^n_+)$. By reordering  the variables, write $x^*=(x_1^*,\dots,x_p^*,0,\dots,0)$   with $x_i^*>0$ for $1\le i\le p$, where $1<p<n$.  Here, the  attractivity of $(x^*,u^*)$ means  the extinction of the populations $x_i(t),\ p+1\le j\le n$, while the first $p$ populations $x_i(t)$ stabilize with time at the 'saturated' value $x_i^*$. For this situation, next result improves Theorem 3.2. Its statement includes Theorems 3.2 and 4.1 as particular cases.

 \proclaim{Theorem 4.2}.  Assume   that $M$ is a P-matrix, let $(x^*,u^*)$ be the saturated equilibrium of (1.1), and suppose that $x^*=(x_1^*,\dots,x_p^*,0,\dots,0)\ (0\le p\le n)$. 
 Write $n_1=p,n_2=n-p$, and the matrices $A=[a_{ij}], |A|=\big [|a_{ij}|\big], A^-=[a_{ij}^-]$   in the form 
 $$A= \left [\matrix {A_{11}&A_{12}\cr A_{21}&A_{22}}\right],\q
 |A|= \left [\matrix {|A_{11}|&|A_{12}|\cr |A_{21}|&|A_{22}|}\right],\q
 A^-= \left [\matrix {A_{11}^-&A_{12}^-\cr A_{21}^-&A_{22}^-}\right],\q
 $$
 where $A_{kl}, |A_{kl}|, A_{kl}^-$ are $n_k\times n_l$ matrices 
 for $k,l=1,2$. Define also
 $$\hat{\cal M}_{11}=diag \Big (\mu_1- c_1{{d_1}\over{e_1}},\dots,\mu_p-c_p{{d_p}\over{e_p}}\Big )-|A_{11}|,\q {\cal M}_{22}^-=diag(\mu_{p+1},\dots,\mu_n)-A_{22}^-.$$
If the matrix
 $$\hat{\cal M}:= \left [\matrix {\hat{\cal M}_{11}&-|A_{12}|\cr -|A_{21}|&{\cal M}_{22}^-}\right],
 \eq(4.2)$$
 is an M-matrix, then $(x^*,u^*)$ is a global attractor for the solutions $(x(t),u(t))$ of (1.1)-(1.2).

  {\it Proof}. 
   The cases $p=n$ and $p=0$ were treated in Theorems 3.2 and 4.1, respectively. Now, consider $0<p<n$. Again, the proof follows along the lines of the proof of Theorem 3.2, so some details are omitted. 
 
 Assume that $\hat{\cal M}$ is an M-matrix. 
  Choose an arbitrarily small $\de_0>0$. Since $\de_0 I_n+\hat {\cal M}$ is a non-singular M-matrix, there is a positive vector $\eta$  such that  
   $(\de_0 I_n+\hat {\cal M})\eta>0$. After a scaling $x_i\mapsto \bar x_i=\eta_i^{-1}x_i, u_i\mapsto \bar u_i=\eta_i^{-1}u_i,\ 1\le i\le n$,  and dropping the bars for the sake of simplicity, we may suppose that $\eta=(1,\dots,1)$. Choose $\de >0$ small and $g$ satisfying (g1)-(g3) and (2.2), with
$$\eqalignno{
\de_0+ \mu_i-c_i{{d_i}\over{e_i}}>(1+\de)\sum_{j=1}^{n} |a_{ij}|,&\q 1\le i\le n_1,&(4.3)\cr
\de_0+\mu_i>(1+\de)\left(\sum_{j=1}^{n_1} |a_{ij}|+\sum_{j=n_1+1}^n a_{ij}^- \right),&\q n_1+1\le i\le n.&(4.4)\cr}
$$

For $i>n_1$, let $\al_i\ge 0$ be such that $b_i+\al_i= \sum_{j=1}^{n_1} a_{ij}x_j^*$, and define $\al_i=0$ for $1\le i\le n_1$.
We now effect the changes $y_i(t)=x_i(t)-x_i^*, v_i(t)=u_i(t)-u_i^*$  for $1\le i\le n$, so we  keep $y_i(t)=x_i(t),v_i(t)=u_i(t)$ for $n_1+1\le i\le n$. Besides (3.14), we obtain
$$y_i'(t)=-(y_i(t)+x_i^*)H_i(t),\q 1\le i\le n,$$
where now
$$
H_i(t)=\al_i+\mu_iy_i(t)+\sum_{j=1}^{n} a_{ij} \int_0^{\infty}\!\!\!K_{ij}(s)y_j(t-s)\,ds
+c_i \int_0^\infty G_i(s)v_i(t-s)\, ds,\q 1\le i\le n.
\eq (4.5)$$

Define $l_i,L_i$ as in (3.15), and recall that  $0\le -l_i\le L_i$  for $i>n_1$. Set $\displaystyle l=\max_{1\le i\le n_1} l_i,\ L=\max_{1\le i\le n} L_i$. We need to prove that $U:=\max (l,L)=0$.

For any $\vare>0$ small, if $t>0$ is sufficiently large  we have
$$\left | a_{ij} \int_0^{\infty}\!\!\!K_{ij}(s)y_j(t-s)\,ds \right |\le |a_{ij} | (\max (l_j,L_j)+\vare),\q 1\le i,j\le n,$$
$$
\int_0^{\infty}\!\!\!K_{ij}(s)y_j(t-s)\,ds\ge -a_{ij}^-(L_j+\vare),\q  n_1+1\le j\le n.$$

Suppose that $U>0$.
If $U=L_i$ or $U=l_i$ for some $i\in \{1,\dots, n_1\}$,  we choose a sequence $t_k\to\infty$ with $y'_i(t_k)\to 0, y_i(t_k)\to L_i$, respectively $y_i(t_k)\to -l_i$, and $(y_{t_k},v_{t_k})\to(\phi,\psi)\in BC\subset UC_g$ as $k\to\infty$. If $U=l_i>0$ for some $i\in \{1,\dots,n_1\}$ and  $y_i(t)$ is eventually monotone, we proceed as in the proof of Theorem 3.2 and easily get a contradiction. Otherwise, $(t_k)$ may be chosen so that $H_i(t_k)=0$. 
We argue as in the proof of Theorem 3.2, and  obtain the estimates (3.21), (3.23), respectively (where now we suppose that $\eta_j=1$ for all $j$).
As in (3.18), we obtain $\nu:=\lim_{k\to\infty} \int_0^\infty G_i(s)v_i(t_k-s)\, ds>-d_iL_i/e_i$ if $y_i(t_k)\to L_i$, and $\nu<d_il_i/e_i$ if $y_i(t_k)\to -l_i$, so we may suppose that $\de_0>0$ in (4.3)-(4.4) was chosen so that $c_i\nu/L_i>-c_id_i/e_i+\de_0$, respectively $c_i\nu/l_i<c_id_i/e_i-\de_0$. 
By taking limits $k\to\infty, \vare_0\to 0^+$, we derive
$$0\ge  \Big (\mu_i-(1+\de)\sum_{j=1}^{n} |a_{ij}|\Big )U+c_i\nu,$$ a contradiction with (4.3).

If $U=L_i$ for some $i\in \{ n_1+1,\dots,n\}$,  we choose a sequence $t_k\to\infty$ with $y_i(t_k)\to L_i, y'_i(t_k)\to 0$,  proceed as above,  and obtain
$$0\ge  \Big (\al_i+\mu_i-(1+\de)\sum_{j=1}^{n_1} |a_{ij}|- (1+\de)\sum_{j=n_1+1}^n a_{ij}^-\Big )U+c_i\nu,\eq(4.6)$$
where now $0\le \phi_i(-s)\le U, \phi_i(0)=U>0$, thus $d_i\int_0^{\infty} e^{-e_is} \phi_i(-s)\, ds>0$,  which implies that $\nu:=\lim_{k\to\infty} \int_0^\infty G_i(s)v_i(t_k-s)\, ds>0$. Therefore, the above estimate contradicts (4.4).\ter
\smal

For equilibria on the boundary of $\R^n_+$, and depending on the sizes and signs of  coefficients $b_i'$s, one might be able to slightly improve the conditions in Theorem 4.2.

 \proclaim{Theorem 4.3}.  Assume  that $M$ is a P-matrix, let $(x^*,u^*)$ be the saturated equilibrium of (1.1), and suppose that $x^*=(x_1^*,\dots,x_p^*,0,\dots,0)\ (1\le p<n)$. Besides the  notations in the statement of Theorem 4.2, we further denote ${\cal A}_{21}=[\tilde a_{ij}]$, where
 $$\tilde a_{ij}=\system {&a_{ij}^-\, ,&\ {\rm if}\ b_i+\sum_{j=1}^{p} a_{ij}^-x_j^*\le 0
 \cr &|a_{ij}|\, ,&\ {\rm if}\ \ b_i+\sum_{j=1}^{p} a_{ij}^-x_j^*>0\cr },\q i=p+1,\dots,n,\, j=1,\dots,p.
 \eq(4.7)$$
If 
 $$\hat{\cal M}:= \left [\matrix {\hat{\cal M}_{11}&-|A_{12}|\cr -{\cal A}_{21}&{\cal M}_{22}^-}\right]
 \eq(4.8)$$
 is an M-matrix, then $(x^*,u^*)$ is a global attractor for the solutions $(x(t),u(t))$ of (1.1)-(1.2).  \vskip 0cm
 In particular, if $A_{21}\ge 0$, $b_i\le 0$  for $p<i\le n$, and $\hat{\cal M}_{11}$ and ${\cal M}_{22}^-$ are M-matrices,  then $(x^*,u^*)$ is globally attractive.
 
 {\it Proof}. Set $p=n_1,n-p=n_2$, $y_i(t)=x_i(t)-x_i^*, v_i(t)=u_i(t)-u_i^*$ for $1\le i\le n$.
 For each $i>n_1$, 
 the function $H_i(t)$  in (4.5) is  given by 
 $$\eqalign{
H_i(t)=\al_i+\mu_iy_i(t)&+\sum_{j=1}^{n_1} a_{ij} \int_0^{\infty}\!\!\!K_{ij}(s)y_j(t-s)\,ds\cr
&+\sum_{j=n_1+1}^{n} a_{ij} \int_0^{\infty}\!\!\!K_{ij}(s)y_j(t-s)\,ds
+c_i \int_0^\infty G_i(s)v_i(t-s)\, ds\cr
\ge \al_i-\sum_{i=1}^{n_1}a_{ij}^+x_j^*&+\mu_iy_i(t)-\sum_{j=1}^{n} a_{ij}^- \int_0^{\infty}\!\!\!K_{ij}(s)y_j(t-s)\,ds+c_i \int_0^\infty G_i(s)v_i(t-s)\, ds\cr
=- (b_i+\sum_{j=1}^{n_1} a_{ij}^-x_j^*)&+\mu_iy_i(t)-\sum_{j=1}^{n} a_{ij}^- \int_0^{\infty}\!\!\!K_{ij}(s)y_j(t-s)\,ds+c_i \int_0^\infty G_i(s)v_i(t-s)\, ds.\cr}
\eq(4.9)$$
For each $i>n_1$, we can use the arguments in the above proof,  with formula (4.6) replaced by the above estimate   if
$(b_i+\sum_{j=1}^{n_1} a_{ij}^-x_j^*)\le 0.$

Now, if $A_{21}\ge 0$ and $b_i\le 0$  for $p<i\le n$,  the matrix in (4.8) becomes
$\hat{\cal M}= \left [\matrix {\hat{\cal M}_{11}&-|A_{12}|\cr0&{\cal M}_{22}^-}\right],$
 which is an M-matrix if and only if $\hat{\cal M}_{11}$ and ${\cal M}_{22}^-$ are M-matrices. 
\ter

\smal

 In applications,  the following corollary is also useful.
 
 \proclaim{Corollary 4.1}. Assume  that $M$  is a P-matrix, let $(x^*,u^*)$ be the unique saturated equilibrium  of (1.1), and $h_i:[0,\infty)\to\R$  continuous functions with $h_i(t)\to 0$ as $t\to\infty\ (1\le i\le n)$. Under the assumptions of Theorems 3.2, 4.2 or 4.3, then 
 all solutions $(x(t),u(t))$ of 
 $$\system{
\displaystyle x_i^{\prime}(t)&=x_i(t)\biggl(b_i-\mu_ix_i(t)-\sum_{j=1}^n a_{ij} \int_0^{\infty}K_{ij}(s)x_j(t-s)\, ds-c_i \int_0^\infty G_i(s)u_i(t-s)\, ds-h_i(t)\biggr), \cr
\displaystyle u_i^{\prime}(t)&=-e_i u_i(t)+d_i x_i(t), \quad i=1,2,\ldots,n , \cr}
\eq(4.10)$$
with initial conditions (1.2) satisfy
 $(x(t),u(t))\to (x^*,u^*)$ as $t\to\infty.$
 
 {\it Proof}. The result follows by repeating  the above proofs with $H_i(t)$ in (3.19), (4.5) or (4.9) replaced by ${\cal H}_i(t):=H_i(t)+h_i(t),\  i=1,\dots,n.$\ter
 \med
\noindent{{\bf Example 4.1}.  We introduce a delayed control in  the single population  model proposed by Volterra  and studied by Miller  [17]:  
$$ \system{
x'(t)&=x(t)\Big ( a-bx(t)-\int_c^t f(t-s) x(s)\, ds-\int_c^t g(t-s) u(s)\, ds\Big ),\cr
u'(t)&=-eu(t)+dx(t),\cr}
\eq(4.11)$$
where $c=0$ or $c=-\infty,\ a,b,d,e>0$, and  the  memory functions  $f:[0,\infty)\to\R, g:[0,\infty)\to [0,\infty)$ are continuous and  in $L^1[0,\infty)$, $g(0)>0$. For $c=-\infty$,  (4.11) is the autonomous system
$$ \system{
x'(t)&=x(t)\Big ( a-bx(t)-\int_0^{\infty} f(s)x(t-s)\, ds-\int_0^{\infty} g(s)u(t-s)\, ds\Big ),\cr
u'(t)&=-eu(t)+dx(t),\cr}
$$
whereas for $c=0$ (4.11) takes the form
$$ \system{
x'(t)&=x(t)\Big ( a-bx(t)-\int_0^{t} f(s)x(t-s)\, ds-\int_0^{t} g(s)u(t-s)\, ds\Big ),\cr
u'(t)&=-eu(t)+dx(t).\cr}
$$
From Theorem  3.2 (see also Remark 3.3) and Corollary 4.1,  if
$b\ge (d/e) \int_0^\infty g(s)\, ds+ \int_0^\infty |f(s)|\, ds$,
then for  any positive solution $(x(t),u(t))$ of (4.11)  with either $c=0$ or $c=-\infty$, we have
$x(t)\to x^*=a\big [b+(d/e) \int_0^\infty g(s)\, ds+ \int_0^\infty f(s)\, ds\big ]^{-1}$ as $t\to\infty$. 

 \med
 When predator-prey systems (1.1) are considered, next result provides less restrictive  sufficient conditions for the extinction of all the predator populations.
   
  \proclaim{Theorem 4.4}.  Assume   that $M$ is a P-matrix, let $(x^*,u^*)$ be the saturated equilibrium of (1.1), and suppose that $x^*=(x_1^*,\dots,x_p^*,0,\dots,0)\ (1\le p< n)$.  With the  notations of Theorem 4.2, assume that $A_{12}\ge 0, A_{21}\le 0$. If $\hat{\cal M}_{11},{\cal M}_{22}^-$ are M-matrices, then
  $(x^*,u^*)$ is a global attractor for the solutions $(x(t),u(t))$ of (1.1)-(1.2). 
  
   {\it Proof.}  Write $n_1=p,n_2=n-p$. Let $(x(t),u(t))$ be a positive solution of (1.1), and set  $y_i(t)=x_i(t)-x_i^*, v_i(t)=u_i(t)-u_i^*$  for $1\le i\le n_1$, and $y_i(t)=x_i(t),v_i(t)=u_i(t)$ for $n_1+1\le i\le n$.  
  \med
  
  {\it Claim 1}. $\displaystyle\limsup_{t\to\infty}x_i(t)\le x_i^*$ for $i=1,\dots,n_1.$
  \med
  
 With $A_{12}\ge 0$,   together with equations (3.14) we get
$$\eqalign{
y_i^{\prime}(t)&\le -(y_i(t)+x_i^*)\biggl(\mu_iy_i(t)-\sum_{j=1}^{n_1} |a_{ij}| \int_0^{\infty}K_{ij}(s)y_j(t-s)\, ds+c_i \int_0^{\infty}G_i(s)v_i(t-s)\, ds\biggr),\cr
 }
 $$
 for $ i=1,2,\ldots,n_1$.
  Fix any $\de_0>0$ small.   With  $\hat{\cal M}_{11}$  an M-matrix, and after a scaling of the variables, we may suppose that $(\de_0I_{n_1}+\hat{\cal M}_{11})\eta>0$ for the positive vector $\eta=(1,\dots,1)\in\R^{n_1}$.  Define $L_i=\limsup_{t\to\infty} y_i(t)$, $U=\max_ {1\le i\le n_1}L_i$. We need to prove that $U\le 0$.  

Suppose that $U>0$.
As for the estimates (3.20), for any $\vare>0$ the definition of $U$ implies that $ \int_0^{\infty}K_{ij}(s)y_j(t-s)\, ds\le (U+\vare)$ for $t>0$ large and $j=1,\dots,n_1$.
Applying the proof of Theorem 3.2, it is clear that  we shall get a contradiction, as in (3.22).
\med

 {\it Claim 2}. $\displaystyle\limsup_{t\to\infty} \int_0^{\infty}K_{ij}(s)y_j(t-s)\, ds\le 0$ for $j=1,\dots,n_1, i=1,\dots,n.$
 
 \med 
 Fix $j\in\{ 1,\dots,n_1\}, i\in\{ 1,\dots,n\}$, and $\de >0$. Since $y_j(t)$ is uniformly bounded in $\R$, there is $T_1>0$ such that  $ \int_{T_1}^{\infty}K_{ij}(s)|y_j(t-s)|\, ds\le \de/2$.
 From Claim 1, $\limsup_{t\to\infty}y_j(t)\le 0$, hence there is $T_2\ge T_1$ such that $y_j(t)<\de/2$ for each $t\ge T_2$. Thus, for $t\ge 2T_2$, we have 
 $$\int_0^{\infty}K_{ij}(s)y_j(t-s)\, ds\le \int_0^{T_2}K_{ij}(s)y_j(t-s)\, ds+\de/2<\de.$$
This proves Claim 2.

\med
  
  {\it Claim 3}. $\lim\limits_{t\to\infty}x_i(t)=0$ for $i=n_1+1,\dots,n.$
  \smal
  
   For each $i\in\{ n_1+1,\dots,n\}$, we only need to prove that $\limsup_{t\to\infty}x_i(t)\le 0$.  Together with the equations $ u_i'(t)=-e_i x_i(t)+d_i u_i(t)$,  we now obtain
 $$\eqalign{
 x_i^{\prime}(t)&=x_i(t)\biggl(b_i-\mu_ix_i(t)-\sum_{j=1}^n a_{ij} \int_0^{\infty}K_{ij}(s)x_j(t-s)\, ds-c_i \int_0^{\infty}G_i(s)u_i(t-s)\, ds\biggr)\cr
 &\le x_i(t)\biggl(\be_i-\mu_ix_i(t)-\!\!\!\sum_{j=n_1+1}^n a_{ij} \int_0^{\infty}K_{ij}(s)x_j(t-s)\, ds-c_i \int_0^{\infty}G_i(s)u_i(t-s)\, ds-h_i(t)\biggr),\cr}$$
 where $\be_i:=b_i-\sum_{j=1}^{n_1}a_{ij}x_j^*\le 0$ (by the definition of a saturated equilibrium) and 
 $$h_i(t)=\sum_{j=1}^{n_1}a_{ij} \int_0^{\infty}K_{ij}(s)y_j(t-s),\q i=n_1+1,\dots,n.$$
From Claim 2 and since $A_{21}\le 0$, we have  $\limsup_{t\to\infty}(-h_i(t))\le 0$ for $i=n_1+1,\dots,n$. From Corollary 4.1 (see also the proof of Theorem 4.2),  the hypothesis ${\cal M}_{22}^-$ is an M-matrix implies Claim 3.

 \smal
  
  {\it Claim 4}. $\lim\limits_{t\to\infty}x_i(t)= x_i^*, \lim\limits_{t\to\infty}u_i(t)= u_i^*$ for $i=1,\dots,n_1.$
  \smal
  
  We write
   $$\system{
y_i^{\prime}(t)&= -(y_i(t)+x_i^*)\biggl(\mu_iy_i(t)+\sum_{j=1}^{n_1} a_{ij}\int_0^{\infty}K_{ij}(s)y_j(t-s)\, ds\cr
 &\hskip 3cm+c_i \int_0^{\infty}G_i(s)v_i(t-s)\, ds+h_i(t)\biggr), \cr
 v_i^{\prime}(t)&=-e_i v_i(t)+d_i y_i(t), \quad i=1,2,\ldots,n_1 , \cr}
$$
where now
$$h_i(t)=\sum_{j=n_1+1}^{n}a_{ij} \int_0^{\infty}K_{ij}(s)x_j(t-s)\, ds,\q i=1,\dots,n_1.$$
Using arguments as the ones above to prove Claim 2, where now we use Claim 3 instead of Claim 1, we get $\lim_{t\to\infty}h_i(t)= 0,\ 1\le i\le n_1$. Claim 4 follows again from Corollary 4.1.\ter

\med

It is straightforward to apply the above results  to uncontrolled systems (3.1), which in the case of saturated equilibria on $\p(\R^n_+)$  lead to better criteria than the ones in [4], as stated below.

\proclaim{Corollary 4.2}.  Assume that $M_0$ is a P-matrix, let $x^*$ be the saturated equilibrium of (3.1), and suppose that $x^*=(x_1^*,\dots,x_p^*,0,\dots,0)\ (1\le p< n)$. With the  notations in the statement of Theorems 4.2, denote also
$$\hat{\cal M}_0:= \left [\matrix {\hat{\cal M}_{0,11}&-|A_{12}|\cr -{\cal A}_{21}&{\cal M}_{0,22}^-}\right],
 $$
where $$\hat{\cal M}_{0,11}=
diag \Big (\mu_1,\dots,\mu_p\Big )-|A_{11}|,\q {\cal M}_{0,22}^-=diag(\mu_{p+1},\dots,\mu_n)-A_{22}^-$$
 and  ${\cal A}_{21}=[\tilde a_{ij}]$ is given by (4.7). If  $\hat{\cal M}_0$  is a non-singular M-matrix, then  $x^*$ is a global attractor for all positive  solutions $x(t)$ of (3.1). \vskip 0cm
Moreover, if either (i) $A_{21}\ge 0$, $b_i\le 0$ for $p<i\le n$, or (ii) $A_{12}\ge 0, A_{21}\le 0$, and $\hat{\cal M}_{0,11},{\cal M}_{0,22}^-$  are non-singular M-matrices, then $x^*$ is a global attractor for the positive  solutions  of (3.1). 

\bigskip

{\bf 5. The two-species Lotka-Volterra system}

\med

As an application of the results in the previous sections, we now analyse with some attention the dynamics for a planar controlled Lotka-Volterra system with delays, without any special constraints on the signs of the Malthusian coefficients $b_i$ and  intra- and inter-specific coefficients $a_{ij}$. For the sake of simplicity, we consider a planar system (1.1) with discrete delays, but the analysis below can be performed for infinite distributed delays as well. 

Consider the system
$$\system{
\displaystyle x_1^{\prime}(t)&=x_1(t)\Bigl(b_1-\mu_1x_1(t)- a_{11} x_1(t-\tau_{11})- a_{12} x_2(t-\tau_{12})-c_1^0 u_1(t)-c_1^1 u_1(t-\sigma_1)\Bigr) \cr
\displaystyle u_1^{\prime}(t)&=-e_1 u_1(t)+d_1x_1(t)\cr
\displaystyle x_2^{\prime}(t)&=x_2(t)\Bigl(b_2-\mu_2x_2(t)- a_{21} x_1(t-\tau_{21})- a_{22} x_2(t-\tau_{22})-c_2^0 u_2(t)-c_2^1  u_2(t-\sigma_2)\Bigr) \cr
\displaystyle u_2^{\prime}(t)&=-e_2 u_2(t)+d_2 x_2(t),\cr}
\eq(5.1)$$
where: $ \mu_i,  c_i^0,  d_i,  e_i$ are positive constants, $c_i^1\ge 0,b_i,a_{ij}\in\R,\tau_{ij}, \sigma_i\ge 0$, $i,j=1,2$.  Denote $c_i=c_i^0+c_i^1,i=1,2$. With the above notation, the community matrix is
$$M=\left [\matrix {\la_1+a_{11}&a_{12}\cr a_{21}&\la_2+a_{22}\cr}\right ] \q  {\rm where }\q \la_i=\mu_i+{{c_id_i}\over {e_i}},\ i=1,2.$$
In what follows, we suppose in addition that $M$ is P-matrix,
i.e.,
$$\det M >0\q {\rm and}\q \la_i+a_{ii}>0,\, i=1,2.
\eq(5.2)$$

There are three possible equilibria on the boundary of $\R^4_+$: the trivial equilibrium $E_0=(0,0,0,0)$, $E_1=({b_1\over {\la_1+a_{11}}}, {{b_1d_1}\over {(\la_1+a_{11})e_1}}, 0,0)$ if $b_1>0$,  and $E_2=( 0,0, {b_2\over {\la_2+a_{22}}}, {{b_2d_2}\over {(\la_2+a_{22})e_2}})$ if $b_2>0$. There is a positive equilibrium
$E^*= (x_1^*,u_1^*,x_2^*,u_2^*)$, where
$$x_1^*={{b_1(\la_2+a_{22})-a_{12}b_2}\over {\det M}},\q  x_2^*={{b_2(\la_1+a_{11})-a_{21}b_1}\over {\det M}},\q u_i^*={{d_i}\over {e_1}}x_i^*,\ i=1,2,$$
if and only if 
$$ b_1(\la_2+a_{22})>a_{12}b_2,\q 
 b_2(\la_1+a_{11})>a_{21}b_1.\eq(5.3)
$$

As already observed, the trivial equilibrium is saturated if and only if  $b_1,b_2\le 0$. In this case, 0 is globally attractive if $\mu_i-a_{ii}^-\ge 0\ (i=1,2)$ and $(\mu_1-a_{11}^-)(\mu_2-a_{22}^-)\ge a_{12}^-a_{21}^-$.
If  $b_i>0$, then $E_i$ is an equilibrium on the boundary of the positive cone, $i=1,2$.


Next, we  give a detailed analysis of the absolute stability, and lack of  it, for the case of $b_1,b_2$ positive. 
The case $b_2\le 0<b_1$ will be studied afterwards.

Let $b_1,b_2>0$, so that the equilibria $E_0,E_1,E_2$ always exist, with $E_0$  unstable. At least one of the conditions in (5.3) is satisfied;  otherwise, we  get $a_{12}b_2\ge b_1(\la_2+a_{22})>0$ and $a_{21}b_1\ge b_2(\la_1+a_{11})>0$ and
 $[(\la_1+a_{11})(\la_2+a_{22})-a_{12}a_{21}]b_1b_2= b_1b_2\det M\le 0$, which is not possible.
 
We now study the stability of $E_1$. Clearly,  a similar analysis can be performed for $E_2$. The characteristic equation for the linearised equation about $E_1=(X_1, {{d_1}\over {e_1}}X_1,0,0)$ is given by
 $$\det \Delta(\la)=0\q {\rm for}\q \Delta(\la)=\la I_4+\left[\matrix {N(\la)&E(\la)\cr 0& C\cr}\right]
 \eq(5.4)$$
 ($I_n$ is the $n\times n$ identity matrix),
 where $X_1={b_1\over {\la_1+a_{11}}}$ and
 $$\displaylines{
 N(\la)=\left[\matrix {X_1(\mu_1+a_{11}e^{-\la \tau_{11}})&X_1(c_1^0+c_1^1e^{-\la \sigma_1})\cr -d_1&e_1\cr}\right],\q
 E(\la)= \left[\matrix{X_1a_{12}e^{-\la \tau_{12}}&0\cr 0&0\cr}\right],\cr
  C= \left[\matrix{-(b_2-a_{21}X_1)&0\cr -d_2&e_2\cr}\right].\cr}$$
 If $b_2(\la_1+a_{11})>a_{21}b_1$, then $b_2-a_{21}X_1>0$, and  $E_1$ is unstable. If $b_2(\la_1+a_{11})\le a_{21}b_1$, 
 the matrix
  $-C$ is stable, therefore $E_1$ is the unique saturated equilibrium. In fact, in this situation, there is no positive equilibrium, but, as  already observed, the condition $b_1(\la_2+a_{22})>a_{12}b_2$ must hold, and  from a dual analysis we would conclude that $E_2$ is unstable. 
  
    When $E_1$ is the unique saturated equilibrium,   conditions (5.2) are not however sufficient to conclude that $E_1$ is a global attractor of all positive solutions for all sizes of the delays $\tau_{11},\sigma_1$. In fact, the  characteristic roots of (5.4) are $\la=-e_2<0, \la =b_2-a_{21}X_1\le 0$ and the solutions of $h(\la)=0$, where
  $$h(\la)=P(\la)+e^{-\la \tau_{11}}Q(\la)+X_1d_1c_1^1e^{-\la \sigma_1}
  \eq(5.5)$$
  with $P(\la)=\la^2+\la (e_1+X_1\mu_1)+X_1(\mu_1e_1+d_1c_1^0), Q(\la)=a_{11}X_1(\la +e_1)$.
  The equation $h(\la)=0$ is the characteristic equation for the system
  $$\system{
  x'(t)=&-X_1\Big [\mu_1x(t)+a_{11}x(t-\tau_{11})+c_1^0u(t)+c_1^1u(t-\sigma_1)\Big ]\cr
  u'(t)=&-[e_1u(t)-d_1x(t)].\cr}
  \eq(5.6)$$
 With $\tau_{11}, \sigma_1=0$, the solutions $\la$ of $h(\la)$ are the eigenvalues of the matrix
  $$-N(0)=-\left[\matrix {X_1(\mu_1+a_{11})&X_1c_1\cr -d_1&e_1\cr}\right],$$
  which has $\det (-N(0))=X_1(\la_1+a_{11})e_1>0$ (from (5.2)) and trace $T_0:=-X_1(\mu_1+a_{11})-e_1$. If $T_0\le0$, then  $E_1$ is stable as an equilibrium of system (5.6) with $\tau_{11}, \sigma_1=0$, otherwise, $E_1$ is unstable. It is particularly difficult to study a second order characteristic equation with two delays, as  equation (5.5), cf.~e.g.~[2,~25] and references therein.
 For instance, fixing $\sigma_1=0$,  if $T_0<0$, in general there is some $\tau^*\ge 0$ such that  the above system (5.6)  is stable for delays $\tau_{11}<\tau^*$, and unstable if $\tau_{11}>\tau^*$, or the stability can change a finite number of times as $\tau_{11}$ increases, and eventually it becomes unstable  -- and therefore, although saturated, $E_1$ becomes unstable also for (5.1). Now assume that
$\mu_1-|a_{11}|-c_1{{d_1}\over {e_1}}\ge 0.$
Then, the trace $T_0$ of $-N(0)$ is always negative, hence $E_1$ is asymptotically stable for system (5.6) with $\tau_{11},  \sigma_1=0$.
Moreover, the matrix
$$\hat N(0)= \left [\matrix {X_1(\mu_1-|a_{11}|)&-X_1c_1\cr -d_1&e_1\cr}\right ]$$
has $\det \hat N(0)=X_1e_1(\mu_1-|a_{11}|-c_1{{d_1}\over {e_1}})\ge 0$ and trace $\hat T_0=X_1(\mu_1-|a_{11}|)+e_1>0$, hence $\hat N(0)$
is an M-matrix [7]. By [5],  it follows that system (5.6) is exponentially stable for all  delays $\tau_{11},\sigma_1>0$.  By Theorem 4.2,  $E_1$ is    the global attractor    of all positive solutions of (5.1) if 
$$\hat {\cal M}= \left [\matrix {\mu_1-|a_{11}|-c_1{{d_1}\over {e_1}}&-|a_{12}|\cr -a_{21}&\mu_2-a_{22}^-\cr}\right ]
$$ is an M-matrix, or, in other words
$$\mu_1-|a_{11}|-c_1{{d_1}\over {e_1}}\ge 0,\ \mu_2-a_{22}^-\ge 0,\
\Big (\mu_1-|a_{11}|-c_1{{d_1}\over {e_1}}\Big)(\mu_2-a_{22}^-)\ge |a_{12}| a_{21}.
\eq(5.7)$$

Assume now (5.3), so that the positive equilibrium $E^*$ exists. For the linearised equation about $E^*$, written as
$$X'(t)=-[DX(t)+L(X_t)],$$
where $D=diag\, (x_1^*\mu_1,e_1,x_2^*\mu_2,e_2)$,
the characteristic equation is given by $\det \Delta(\la):=\la I_4+D+L(e^{\la \cdot}I_4)=0$, and
similar computations as the ones above lead to 
$$D+L(e^{\la \cdot}I_4)=\left[\matrix {N_1(\la)&E_1(\la)\cr E_2(\la)& N_2(\la)\cr}\right],
 \eq(5.8)$$
 where
$$N_i(\la)=\left[\matrix {x_i^*(\mu_i+a_{ii}e^{-\la \tau_{ii}})&x_i^*(c_i^0+c_i^1e^{-\la\sigma_i})\cr -d_i&e_i\cr}\right],\
 E_i(\la)= \left[\matrix{x_i^*a_{ij}e^{-\la \tau_{ij}}&0\cr 0&0\cr}\right],\q i,j=1,2, j\ne i.$$
 One can easily check that
$\det\Delta(0)=\det ( D+L(I_4))=x_1^*x_2^*e_1e_2\det M$, thus $\det\Delta(0)>0$ since $M$ is a P-matrix. As for the study of the stability of $E_1$, even if $E^*$ is asymptotically stable for the corresponding ODE system obtained by taking all the delays equal to zero in (5.1),  the positive equilibrium  $E^*$ of (5.1) might become unstable as  the delays increase. In fact, by letting $c_1,c_2\to 0^+$, from (5.8) we obtain
$\det \Delta(\la)\to (\la+e_1)(\la+e_2)h(\la)$,
where now
$$h(\la)= \left|\matrix {\la+x_1^*(\mu_1+a_{11}e^{-\la \tau_{11}})&x_1^*a_{12}e^{-\la \tau_{12}}\cr
x_2^*a_{21}e^{-\la \tau_{21}}&\la+x_2^*(\mu_2+a_{22}e^{-\la \tau_{22}})
}\right|.$$
Choosing e.g. $\tau_{ii}=0,(i=1,2)$ and $a_{12}=1,a_{21}=-1$, one can see that it is possible to choose the other coefficients in such a way   that $(\mu_1+a_{11})x_1^*=(\mu_2+a_{22})x_2^*=:b$ and $x_1^*x_2^*=:c>b^2$. Then,
$h(\la)=(\la+b)^2+ce^{-\la (\tau_{12}+\tau_{21})}$,  which has roots $\pm i\sqrt {c-b^2}$ if $\tau:=\tau_{12}+\tau_{21}=\tau_n$, where  
$$\tau_n\in (0,\pi)+2n\pi,\q \tan  ( \tau_n\sqrt {c-b^2})={{2b\sqrt {c-b^2}}\over {c-2b^2}} ,\q  n=0,1,2,\dots.$$
In particular, for $\tau>\tau_0$ and close to $\tau_0$, there is a pair of characteristic roots with positive real parts, thus the equilibrium becomes unstable. Moreover, system (5.1) has a sequence of Hopf bifurcations at $\tau=\tau_n, n=0,1,2,\dots$ [23].
On reverse, if 
$$\hat M=\left [\matrix {\mu_1-|a_{11}|-c_1{{d_1}\over {e_1}}&-|a_{12}|\cr -|a_{21}|&\mu_2-|a_{22}|-c_2{{d_2}\over {e_2}}\cr}\right ]$$
is an  M-matrix, so that we have the conditions
$$(\mu_1-|a_{11}|-c_1{{d_1}\over {e_1}})(\mu_2-|a_{22}|-c_2{{d_2}\over {e_2}})\ge |a_{12}a_{21}|\ {\rm and}\ \mu_i-|a_{ii}|-c_i{{d_i}\over {e_i}}\ge0,\, i=1,2,
\eq(5.9)$$ 
the positive equilibrium $E^*$ is globally attractive, for all sizes of delays $\tau_{ij},\sigma_i$.

As an application of the use of the controls, in the example below we  change the position of the globally attractive equilibrium, from  the boundary  to the interior of $\R^2_+$, recovering one of the species, otherwise condemned to extinction.

\med

\noindent{\bf Example 5.1}. Consider the following uncontrolled system  with $n=2$ and
e.g.~$b_1=1,b_2={1\over 3}, \mu_1=\mu_2=1, a_{11}=a_{22}= a_{21}={1\over 2}, a_{12}={1\over 8}$:
$$\system{
\displaystyle x_1^{\prime}(t)&=x_1(t)\Bigl(1-x_1(t)- {1\over 2} x_1(t-\tau_{11})-  {1\over 8} x_2(t-\tau_{12})\Bigr) \cr
\displaystyle x_2^{\prime}(t)&=x_2(t)\Bigl({1\over 3}-x_2(t)-{1\over 2} x_1(t-\tau_{21})- {1\over 2} x_2(t-\tau_{22})
\Bigr) \cr
}.
$$
With the above notations, we have
$M_0= \left [\matrix{3/2& 1/8\cr 1/2&3/2\cr}\right ],\ \hat M_0= \left [\matrix{ 1/2& -1/8\cr -1/2&1/2\cr}\right ]$.
Its saturated equilibrium  is $(X_1,0)=({2\over 3},0)$. Furthermore, 
$\det M_0>0$ and $\hat M_0$ is a non-singular M-matrix, hence from [4] we derive that  $(X_1,0)$ is GAS.
We now introduce the controls, in order to recover the $x_2(t)$ population, which otherwise would become extinct with time. Clearly, for any choice of positive coefficients $c_i,d_i,e_i,\,  i=1,2$, conditions (5.3) hold, and therefore the controlled system (5.1) with the above coefficients has a positive equilibrium $E^*$. Now, if we choose e.g. $\al_i:=c_i{{d_i}\over {e_i}}\le {1\over 4},\, i=1,2$, we have that
$$\hat M= \left [\matrix{ 1/2-\al_1\ & -1/8\cr -1/2&1/2-\al_2\cr}\right ]$$
is an M-matrix. Invoking Theorem 3.2, we get that
 $E^*$ is a global attractor of all positive solutions.

\med

Now, suppose that $b_2\le 0<b_1$ and $a_{21}\ge 0$. Clearly (5.3) fails to be true,  $E_1$ is the saturated equilibrium, and by Theorem 4.3  $E_1$ is a global attractor of all positive solutions of (5.1) if $$\mu_1-|a_{11}|-c_1{{d_1}\over {e_1}}\ge 0,\q \mu_2-a_{22}^-\ge 0.\eq(5.10)$$ 

Next, consider a typical predator-prey system (5.1), where $b_2< 0<b_1$ and $a_{12}>0, a_{21}< 0$. In the absence of the positive equilibrium, which amounts to have $b_2(\la_1+a_{11})\le a_{21}b_1$, then $E_1$ is the saturated equilibrium. Using now Theorem 4.4, if (5.10) is satisfied, then again $E_1$ is a global attractor. In this framework, we again illustrate how the controls can be used to change the position of a globally attractive saturated equilibrium.

\med

\noindent{\bf Example 5.2}. For the particular case of 
 $b_1=1,b_2=-{5\over 4}, \mu_1=\mu_2=1, a_{11}=a_{22}={1\over 2},a_{12}= {1\over 8}, a_{21}=-2$,
we obtain the predator-prey system without controls
$$\system{
\displaystyle x_1^{\prime}(t)&=x_1(t)\Bigl(1-x_1(t)- {1\over 2} x_1(t-\tau_{11})-  {1\over 8} x_2(t-\tau_{12})\Bigr) \cr
\displaystyle x_2^{\prime}(t)&=x_2(t)\Bigl(-{5\over 4}-x_2(t)+2 x_1(t-\tau_{21})- {1\over 2} x_2(t-\tau_{22})
\Bigr) \cr
},
$$
with community matrix $M_0= \left [\matrix{3/2& 1/8\cr -2&3/2\cr}\right ]$.
For this system,  $(x_1^*,x_2^*)=({{53}\over {80}}, {1\over {20}})$ is the positive equilibrium.
 Moreover, since $\det M_0>0$ and 
$\hat M_0= \left [\matrix{ 1/2& -1/8\cr -2&1/2\cr}\right ]$
is an { M-matrix}, from [4] it follows that  $(x_1^*,x_2^*)$ is globally atractive.
We now introduce the controls, in order to drive the predators to extinction. For the above chosen coefficients,
  $ b_2(\mu_1+a_{11}+c_1{{d_1}\over {e_1}})\le a_{21}b_1$ if and only if $c_1{{d_1}\over {e_1}}\ge {1\over {10}}$, in which case
 $ E_1=({1\over {{3\over 2}+c_1{{d_1}\over {e_1}}}}, 0, {{d_1e_1}\over {{3\over 2}e_1+c_1d_1}},0)$ is the saturated equilibrium.
If we now choose ${1\over {10}}\le c_1{{d_1}\over {e_1}}\le{1\over 2},$  Theorem 4.4  yields that
 $ \mu_1-|a_{11}|-c_1{{d_1}\over {e_1}}\ge 0$,
thus   $E_1$ is a global attractor of all positive solutions.

\med

We  summarise the above global asymptotic behaviour results as follows:

\proclaim{Proposition 5.1}. Consider system (5.1), and assume (5.2).\vskip 0cm
(i) If $b_1,b_2\le 0$, then 0 is the saturated equilibrium; in this case, 0 is globally attractive if $\mu_i-a_{ii}^-\ge 0\ (i=1,2)$ and $(\mu_1-a_{11}^-)(\mu_2-a_{22}^-)\ge a_{12}^-a_{21}^-$.\vskip 0cm
(ii) If (5.3) holds, there exist a positive equilibrium, which is GAS under the additional conditions (5.9).
 \vskip 0cm
(iii)  If $b_1,b_2>0$ and $b_2(\la_1+a_{11})\le a_{21}b_1$, then $E_1$ is the saturated equilibrium; in this case, $E_1$ is a global attractor of all positive solutions if conditions (5.7) are satisfied.
 \vskip 0cm
(iv) If  $b_2\le 0<b_1$, and: (a) either $a_{21}\ge 0$, or (b) $a_{12}>0, a_{21}< 0, b_2(\la_1+a_{11})\le a_{21}b_1$, then $E_1$ is the saturated equilibrium; in this case, $E_1$ is a global attractor of all positive solutions if conditions (5.10) are satisfied.

\bigskip

\centerline{\bf References} 

\medskip

\baselineskip=13.5pt

\item{1.} A. Berman  and   R. Plemmons.  {\it Nonnegative Matrices in the Mathematical Sciences}(Academic Press, New York, 1979).

\item{2.} S.A. Campbell and  Y. Yuan. Zero singularities of codimension two and three in delay differential equations. {\it Nonlinearity} {\bf 21}  (2008), 2671--2691.

\item{3.} F. Chen. The permanence and global attractivity of Lotka-Volterra competition system with feedback controls. {\it Nonlinear Anal. RWA} {\bf 7}  (2006),
133--143.
 
 \item{4.} T. Faria. Stability and Extinction for Lotka-Volterra Systems with Infinite Delay. {\it J. Dynam. Differential Equations} {\bf 22} (2010), 299--324.

\item{5.} T. Faria  and   J.J. Oliveira. Local and global stability for   Lotka-Volterra systems 
with distributed delays and instantaneous feedbacks.
{\it J. Differential Equations} {\bf 244}   (2008),
 1049--1079.
 
 \item{6.}  T. Faria  and   J.J. Oliveira. General criteria for asymptotic and exponencial stabilities of neural network models with unbounded delay. {\it Appl. Math. Comput.}  {\bf 217}  (2011), 9646--9658.

\item{7.}   M. Fiedler. {\it Special Matrices and Their Applications in Numerical Mathematics} (Martinus Nijhoff
Publ., Kluwer, Dordrecht, 1986).

\item{8.}  K. Gopalsamy and P. Weng. Global attractivity in a competition system with feedback controls. {\it Comput. Math. Appl.} {\bf 45}   (2003), 665--676.

\item{9.}    J. Haddock and W. Hornor. Precompactness and convergence in norm of positive orbits in a certain
fading memory space.  {\it Funkcial. Ekvac.} {\bf 31}  (1988), 349--361.

\item{10.}   J. Haddock, M.N.  Nkashama  and  J.H. Wu. Asymptotic constancy for linear neutral Volterra integrodifferential equations.  {\it Tohoku Math. J.} {\bf 41} (1989), 689--710.

\item{11.}  J.K. Hale. {\it Asymptotic Behavior of Dissipative Systems}
(Amer.~Math.~Soc., Providence, 1988).

\item{12.}  J.K. Hale and  J. Kato. Phase space for retarded equations with infinite delay.  {\it Funkcial.
Ekvac.} {\bf 21}  (1978), 11--41.

\item{13.} Y. Hino, S.  Murakami and  T. Naito. {\it Functional Differential
Equations with Infinite Delay} (Sprin\-ger-Verlag, New-York, 1993).

  \item{14.}  J. Hofbauer and  K. Sigmund. {\it The Theory of Evolution and Dynamical Systems} (London
Mathematical Society, Cambridge University Press, Cambridge, 1988). 


\item{15.}   Y. Kuang. Global stability in delay differential systems without dominating instantaneous negative feedbacks. {\it J. Differential Equations} {\bf 119} (1995), 503--532.

\item{16.}  Z. Li, M.  Han and  F. Chen. Influence of feedback controls on an autonomous Litoka-Volterra competitive system with infinite delays. {\it Nonlinear Anal. RWA} {\bf 14} (2013), 402--413.


\item{17.}   R.K. Miller. On Volterra's population equation.
 {\it SIAM J. Appl. Math.} {\bf 14}   (1966), 446--452.
 
 \item{18.}  F. Montes de Oca and   M. Vivas. Extinction in a two dimensional Lotka-Volterra system with infinite delay.  {\it Nonlinear Anal. RWA} {\bf 7}  (2006), 1042--1047.

\item{19.}   Y. Muroya. Partial survival and extinction of species in nonutonomous Lotka-Volterra systems with delays. {\it Dynamic Systems and Applications}  {\bf 12} (2003), 295--306.

\item{20.}  L. Nie, J. Peng and  Z. Teng. Permanence and stability in multi-species non-autonomous LotkaÐVolterra competitive systems with delays and feedback controls. {\it Math. Comput. Modelling} {\bf 49} (2009), 295--306.


\item{21.} C. Shi, Z. Li  and  F. Chen. Extinction in a nonautonomous Lotka-Volterra competitive system with infinite delay and feedback controls. 
{\it Nonlinear Anal. RWA} {\bf 13} (2012), 2214--2226.

\item{22.}  H.L. Smith. {\it Monotone Dynamical Systems. An
Introduction to the Theory of Competitive and Cooperative Systems}
 (Amer. Math. Soc., Providence, 1995).

\item{23.}  P. T\'aboas. Periodic solutions of a planar delay equation. 
{\it Proc. Roy. Soc. Edinburgh Sect.~A} {\bf 116}  (1990), 85--101.

\item{24.}  P. Weng.  Existence and global stability of positive periodic solution of periodic integrodifferential systems with feedback controls. {\it Comput. Math.
Appl.} {\bf 40}  (2000), 747--759.

\item{25.}  D. Xiao and  S. Ruan. Multiple bifurcations in a delayed predator-prey system with nonmonotonic functional response. {\it J. Differential Equations } {\bf 176}  (2001), 494--510.

\item {26.}   R. Xu, M.A.J. Chaplain and F.A. Davidson. Global asymptotic stability in a nonautonomous
$n$-species Lotka-Volterra predator-prey system with infinite delays. {\it Appl. Anal.} {\bf 80}(2002), 107--126.

\item{27.} Z. Zhang. Existence and global attractivity of a positive periodic solution for a generalized delayed population model with stocking and feedback control. {\it Math. Comput. Modelling} {\bf 48} (2008), 749--760.

\end